\documentclass[3p]{elsarticle}

\usepackage{amsmath}
\usepackage{amssymb}
\usepackage{color}
\usepackage{hyperref}
\usepackage{placeins}

\journal{Journal of Computational Physics}

\begin{document}

\begin{frontmatter}

\title{A Fully Lagrangian Meshfree Framework for PDEs on Evolving Surfaces}

\author[affil1]{Pratik Suchde\corref{mycorrespondingauthor}}
\ead{pratik.suchde@itwm.fraunhofer.de}
\author[affil1]{J\"org Kuhnert}

\cortext[mycorrespondingauthor]{Corresponding author}

\address[affil1]{Fraunhofer ITWM, 67663 Kaiserslautern, Germany}

\begin{abstract}
We propose a novel framework to solve PDEs on moving manifolds, where the evolving surface is represented by a moving point cloud. This has the advantage of avoiding the need to discretize the bulk volume around the surface, while also avoiding the need to have a global mesh. Distortions in the point cloud as a result of the movement are fixed by local adaptation. We first establish a comprehensive Lagrangian framework for arbitrary movement of curves and surfaces given by point clouds. Collision detection algorithms between point cloud surfaces are introduced, which also allow the handling of evolving manifolds with topological changes. We then couple this Lagrangian framework with a meshfree Generalized Finite Difference Method (GFDM) to approximate surface differential operators, which together give a method to solve PDEs on evolving manifolds. The applicability of this method is illustrated with a range of numerical examples, which include advection-diffusion equations with large deformations of the surface, curvature dependent geometric motion, and wave equations on evolving surfaces.
\end{abstract}


\end{frontmatter}


\section{Introduction}

The need to numerically solve PDEs on evolving surfaces arises in various fields. From the modelling of surfactants \cite{Stone1990} and other fluid dynamics \cite{Gilman2000}, to the modelling of airbags \cite{Hirth2007} and parachutes \cite{Accorsi2000}. PDEs on deforming curves and surfaces also appear in the modelling of biomembranes \cite{Elliott2010} and cell motility \cite{Elliott2012}, in visualization \cite{Kim2013}, and image processing \cite{Jin2004}. 

Most methods for solving PDEs on moving surfaces are mesh-based. Both finite elements  \cite{Dziuk2007,Dziuk2013,Lubich2013} and finite volumes \cite{Nemadjieu2014} have been used in this context. The movement of the surface is achieved by moving the discretizing mesh. This introduces the trouble of mesh distortion which needs an expensive remeshing \cite{Li2018}. Note that this is the same issue with moving mesh methods for volumetric flow with free surfaces. To avoid the need of meshing and remeshing, surfaces are often discretized with a cloud of numerical points. As a result, a significant amount of work has been done towards the development of meshfree methods for solving PDEs on static surfaces (for example, \cite{Flyer2014, Liang2013,Marz2012}). Some of these meshfree methods have also been extended to handle evolving surfaces.

To the end of solving PDEs on stationary surfaces, a lot of work has been done to derive meshfree methods that scale with the dimension of the surface itself, rather than the dimension of the embedding space \cite{Fuselier2013,Liang2013}. This significantly decreases the computational cost as compared to methods which discretize the embedding space (or a subset of the same dimension). However, when it comes to PDEs on evolving surfaces, most meshfree methods require a discretization of the bulk volume surrounding the surface \cite{Auer2013,Leung2009,Petras2019, Petras2016,Sokolov2017}. The movement of the discrete manifold is captured by some notion of tracking on this surrounding grid or point cloud. This is similar in essence to interface tracking for free surface volumetric flows done in static mesh-based methods \cite{Gueyffier1999}, and all inaccuracies in that context carry over to the present context. Thus, by requiring tracking of the surface, these methods lose one of the fundamental advantages of meshfree methods. Moreover, they require a discretization of a higher dimensional space, which makes them computationally expensive. 
 
In the context of volumetric flows, it is well known that Lagrangian frameworks can accurately capture advection as well as the shape of interfaces and free surfaces \cite{Lamb1993,Samulyak2018}. In this paper, we introduce such a fully Lagrangian framework to accurately capture the movement of an evolving surface. We present a new meshfree method for solving PDEs on evolving surfaces that scales with the manifold dimension. The surface is discretized with a cloud of numerical points, without any need to discretize the embedding space around the surface. Movement of the surface is captured using a fully Lagrangian framework. Spatial derivatives are computed by virtual projections to the tangent space at each point. A major advantage of a moving point cloud over a moving mesh is that point cloud distortion is much easier to fix. For point clouds defining volumes, this has been shown in, for example, \cite{Drumm2008}. Here, we will show that the same also holds good for surfaces. To this end, we present methods to fix distortions in surface point clouds by adding and removing points where necessary by purely local considerations.

We note that moving Lagrangian particles have already been used to solve PDEs on evolving surfaces by several authors. The novelty in the present work lies in the fact that only the surface is being discretized, and not the bulk around it. In contrast, \cite{Bergdorf2010} uses Lagrangian particles in a band around the manifold, and a regular grid in the background to interpolate particle locations and to fix distortion. While \cite{Leung2011} and \cite{Petras2016} use Lagrangian particles for the moving surface with a regular fixed grid of the dimension of the embedding space in the background for reference and neighbourhood information. 

%

The paper is organized as follows. Section \ref{sec:Discretization} contains some information about the setup of surface point clouds. Section~\ref{sec:LagrangianFramework} describes how a Lagrangian framework can be established for surface PDEs, and includes the details about the Lagrangian movement, the required adaptation of nodes, and contact handling. Section~\ref{sec:PDEsMovingMan} then goes on to show how this can be used to solve PDEs on moving manifolds in a Lagrangian way, and presents numerical examples. The paper is then concluded with a short discussion in Section~\ref{sec:Conclusion}.

\section{Preliminaries}
\label{sec:Discretization}

To distinguish between the cases of PDEs on surfaces and those on volume domains, we use the term `volumetric' to denote the volume domain case. 

Throughout this work, we establish initial point clouds as per \cite{Drumm2008}. Notation and definitions for the point clouds are in accordance with \cite{Suchde2019,Suchde2017_CCC}. Conventions typically used for volumetric meshfree GFDM point clouds \cite{Drumm2008, Jefferies2015} are carried over to the surface case here. Throughout this paper, we consider only smooth orientable $2$ manifolds in $\mathbb{R}^3$. However, the ideas presented in this paper can easily be generalized to higher dimensions and co-dimensions.

The smooth oriented surface or manifold $M$ is discretized with $N$ non-uniformly spaced numerical points also referred to as nodes or particles. These points are simply locations where approximations are carried out, and they do not carry mass. The $N$ points include points both in the interior and at the boundary~(if any) of the manifold. Approximations at a numerical point $i=1,2,\dots,N$ are done based on a support or neighbourhood $S_i$ of $n_{S_i}$ nearby points, within a distance of $h$. $h=h(\vec{x},t)$ is referred to as the interaction radius or smoothing length. The spatial distribution of points is described by three parameters: $h$, $r_{max}$, and $r_{min}$. It is ensured that there is no hole of size $r_{max}h$ within the point cloud, and that no two points are closer than $r_{min}h$. Both $r_{max}$ and $r_{min}$ are global constants, taken to be the same for all simulations, and are not dependent on the PDE being solved. As a result, the smoothing length $h$ serves not just as the size of the support at each point, but also as the spatial discretization size. $r_{max}$ and $r_{min}$ determine the number of neighbouring points in each support domain. We use $r_{max}=0.45$ and $r_{min}=0.2$, which result in about $15-20$ points in each neighbourhood. These values are carried over from conventions on point cloud spacing used in meshfree GFDM volumetric flow solvers \cite{Drumm2008, Jefferies2015, Suchde2017_CCC}. All distance computations are done in the embedding space (here, in $\mathbb{R}^3$), and \textit{not} along the manifold. To satisfy these maximum and minimum distance criteria on manifolds, which could change with time, addition and deletion of points on the surface is carried out, which is explained in Sections~\ref{sec:PCMaddition} and \ref{sec:PCMdeletion}.

\section{The Lagrangian Framework}
\label{sec:LagrangianFramework}

In this section we present details about how a meshfree Lagrangian framework can be established for surfaces. To the best of the authors' knowledge, such a fully Lagrangian meshfree setting without any background grid has never been done for manifolds.

We consider an advection velocity $\vec{v}$, which can have components both normal and tangential to the surface. In contrast to a lot of existing methods for PDEs on evolving surfaces
, we do not split the velocity into its normal and tangential components, and use the entire velocity as the Lagrangian transport velocity. Furthermore, unlike many mesh-based Lagrangian methods for surface PDEs
, we do not assume any correspondence between the point clouds at different time levels. At a given time level, no mapping is preserved to the initial configuration of points, and all approximations are done at the present point cloud directly.

\subsection{The Actual Movement}

The Lagrangian motion step involves solving the ODE system

\begin{equation}
	\label{Eq:LagrangianMotion}
	\frac{d\vec{x}}{dt} = \vec{v} \,,
\end{equation}
where $\vec{v}$ is the advection velocity. In most volumetric Lagrangian frameworks for meshfree methods, Eq.\,\eqref{Eq:LagrangianMotion} is solved by a first order method which assumes the velocity is constant between two time levels. The same is also done in Lagrangian and semi-Lagrangian methods for surfaces. Both in the moving mesh-context \cite{Li2018}, and meshless manifold tracking context \cite{Petras2016}. This involves each point of the point cloud, or equivalently each node of the mesh, being moved with the given velocity field as
\begin{equation}
	\label{Eq:Move_o1}
	\vec{x}^{(n+1)} = \vec{x}^{(n)} + \vec{v}^{\,(n)}\Delta t \,,
\end{equation}
where the bracketed superscripts indicate the time level. Here, time-integration is done between the time levels $t^n$ and $t^{n+1}$, and it is assumed that $\vec{v}^{\,(n+1)}$ is unknown. Alternatively, similar methods are also done by moving with the velocity $\vec{v}^{\,(n+1)}$, if known, or an average of $\vec{v}^{\,(n+1)}$ and $\vec{v}^{\,(n)}$. In each case, the velocity is taken to be constant within each time step. 

In our earlier work for meshfree volumetric Lagrangian flows \cite{Suchde2018_PCM}, we have shown the inaccuracies surrounding the first order movement similar to Eq.\eqref{Eq:Move_o1}, which lead to large defects in volume conservation. To get accurate Lagrangian movement with Eq.\,\eqref{Eq:Move_o1}, most authors tend to very small time steps. The same arguments of inaccuracy presented in \cite{Suchde2018_PCM} for volumetric flows carry over to the surface case here. Thus, we use the more accurate second order method for point cloud movement
\begin{equation}
	\label{Eq:Move_o2}
	\vec{x}^{\,(n+1)} = \vec{x}^{\,(n)} + \vec{v}^{\,(n)}\Delta t + \frac{1}{2}\frac{\vec{v}^{\,(n)} - \vec{v}^{\,(n-1)} }{\Delta t} \Delta t ^2\,.
\end{equation}
Alternatively, if $\vec{v}^{\,(n+1)}$ is known at the time of updating point locations, $\vec{v}^{\,(n)}$ and $\vec{v}^{\,(n-1)}$ can be replaced by $\vec{v}^{\,(n+1)}$ and $\vec{v}^{\,(n)}$ respectively in Eq.\,\eqref{Eq:Move_o2}. This would be the case if the new velocity has been determined before the movement step. Note that since Eq.\,\eqref{Eq:Move_o2} retains the explicit nature of Eq.\,\eqref{Eq:Move_o1}, the second order method comes at virtually no extra cost over the first order one. 

To emphasize the need of such a second order method according to Eq.\,\eqref{Eq:Move_o2}, a simple case of a rotating quarter-sphere is considered in Section~\ref{sec:RotatingQSphere}.

We note that most practical applications for PDEs on evolving surfaces will rely on time-dependent velocity fields. While even higher order integration methods for the movement step can be developed for time-independent velocity fields, they can not be easily generalized to the time-dependent case for moving point clouds~\cite{Suchde2018_PCM}, so they are not very interesting in the present context. 

In both Eq.\,\eqref{Eq:Move_o1} and Eq.\,\eqref{Eq:Move_o2}, each point of the point cloud is moved explicitly and independently of all other points. Here, it is assumed that any implicit or global constraint on the motion of the surface is already represented in $\vec{v}$.

\subsection{Normal Computation}
\label{sec:NormalComputation}
Each interior point $i$ of the manifold must be equipped with a unit normal $\vec{n}_i$ and two unit tangents $\vec{t}_{1,i}$ and $\vec{t}_{2,i}$. Furthermore, each boundary point $i$ (if any) of the manifold must be equipped with a unit manifold normal $\vec{n}_i$, a unit tangent $\vec{t}_i$, and a unit boundary normal $\vec{\nu}_i$. Note that $\vec{\nu}_i$ is normal to the boundary of the manifold, but lies in the tangent space of the manifold. At each point, these normals and tangents form an orthonormal system of vectors.

Multiple possibilities exist for the computations of these normals and tangents on point clouds. A popular way to determine these is by using principal component analysis~(PCA) \cite{Liang2013,Mitra2003} which computes eigenvalues of local covariance matrices. 

PCA is known to be sensitive to outliers and noise in the point cloud. Weighted PCA approaches have been used to overcome the same \cite{Petronetto2013}. 
To use a more robust method to compute normals on point clouds, we use a minimization approach that maximizes the angles between the normal and the neighbouring points. 
This is based on the procedure used for computing normals at the free surface for volumetric flow by \cite{Edgar2017}. For an interior point $i$, the manifold normal $\vec{n}_i$ is set to be the unit normal vector minimizing the functional
\begin{equation}
	\sum_{j\in S_i} \cos^2(\theta_{ij})W_{ij} = \vec{n}_i^{\,T} \sum_{j\in S_i}\left( W_{ij} \frac{ ( \vec{x}_j - \vec{x}_i)( \vec{x}_j - \vec{x}_i)^T }{\| \vec{x}_j - \vec{x}_i\|^2} \right) \vec{n}_i \,,
\end{equation}
where $\theta_{ij}$ is the angle between the normal $\vec{n}_i$ and the distance vector $\vec{x}_j - \vec{x}_i$, and $W_{ij}$ is a Gaussian weighting function to emphasize the effect of the closest neighbours. Once $\vec{n}_i$ is known, the tangents can be computed so as to obtain an orthonormal system of vectors.

For boundary points, the boundary normal $\vec{\nu}_i$ is first computed considering only the boundary neighbours. Based on this, the boundary tangent is computed. After which, the manifold normal can be established at the boundary point.

This process of normal and tangent computation needs to be done at each time step, since the positions of each point are time-dependent.

We emphasize that we consider surfaces that are orientable, and the computed normal field should also result in an oriented surface. For example, for the case of the surface of a sphere, all the normals must either be pointing inwards, or all outwards. The need of an oriented normal field becomes especially important, for example, when computing the curvature of the surface by the surface divergence of the normal field $\kappa = -\frac12\nabla_M\cdot\vec{n}$.

The same is explicitly checked by a sweep through all points when the normals are computed for the first time during a simulation. At later time steps, to maintain the orientation of the normal field, newly computed normals are compared with the old values. The orientation of $\vec{n}_i^{\,\text{new}}$ is chosen such that
\begin{equation}
	\vec{n}_i^{\,\text{new}} \cdot \vec{n}_i^{\,\text{old}} > 0 \,,
\end{equation}
where the superscripts $new$ and $old$ indicate the normals at the current and previous time levels respectively. Similar considerations are also done for normal computation at newly added points~(to be explained in a coming section), but based on the normals of neighbouring points.

\subsection{Neighbour Searching}
\label{sec:NbrSearch}

Efficient neighbour searching is crucial for the efficiency of meshfree methods for surface PDEs, just as it is in the volumetric case. Neighbour search algorithms used in volumetric meshfree methods \cite{Dominguez2010,Drumm2008, Onderik2008} directly carry over to the surface case being considered here. For each node, its neighbouring nodes are determined as the nodes within a certain distance $h$ from it. The naive approach to neighbour searching would involve computing distances between every pair of nodes in the computational domain. However, this procedure is exorbitantly expensive. To avoid excessive distance computations, the domain is usually split into multiple regions referred to as boxes or cells. Using such a decomposition, for each node, distances only need to be computed with other nodes within the same box (or possibly also adjacent boxes). Several different data structures have been used to this end. One of the ways to do the same is to use quadtree or octree type searching algorithms. For the Lagrangian case, it needs to be ensured that  when a new point is added, it is added in the correct box. To save time, it is often desirable to preserve the tree structure for a few time steps, and reconstruct it only every $p$ time steps.


%
%
\subsection{Addition of Points}
\label{sec:PCMaddition}

Just like the mesh-based case, Lagrangian movement of a surface point cloud can cause distortions which manifest as `holes' in the point cloud where no points are present, or `clusters' where too many points are present. To fix this, points are added in holes, and removed or merged when they come too close. We now propose a method for hole identification and filling in the absence of a background discretization of the embedding space.

\begin{figure}
  \centering
  \includegraphics[width=\textwidth]{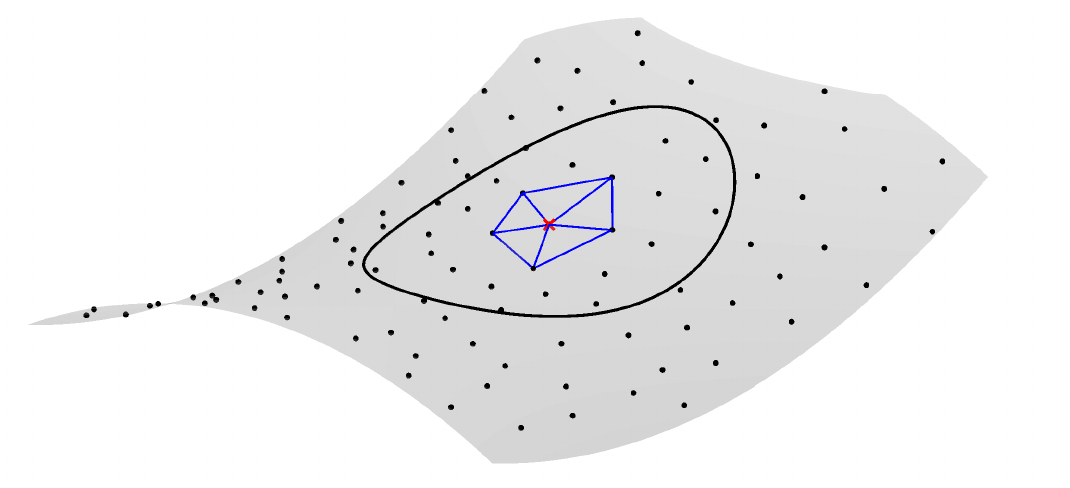}
  \caption{Local triangulation of a point on a surface. The central point is marked with a red cross. The black circle around it marks the neighbourhood of the central point. The blue triangles indicate the triangulation computed to determine locations where points need to be added to maintain regularity of the point cloud. The grey region indicates the surface specified by the black points.}
  \label{Fig:LocalTriangulation}%
\end{figure}

Several volumetric meshree methods use locally defined meshes~\cite{Liu2009Meshfree} for a variety of purposes, from integration in weak form methods~\cite{Atluri1998}, to support domain selection~\cite{Liu2013PIM}, to post-processing~\cite{Jefferies2015}. We propose to use locally defined triangulations to identify holes in the surface point cloud. For this, at each point, we compute a `one-ring' of triangles based solely on the locations of neighbouring points. An example of such a triangulation within a support domain is shown in Figure~\ref{Fig:LocalTriangulation}. 

It is important to note the difference between the above mentioned locally computed triangulation, and the generation of an actual mesh on the entire surface. The local triangulations are much easier to generate. They scale very well when computed in parallel as the procedure for each point is independent of that for all other points. There is no imposed restriction for the triangles to be `good'. All the local 1-ring of triangles need not stitch together to form a global mesh of the surface. In fact, in most cases, they do not, unless the support domains are considered very large, which is rarely done in practice. Furthermore, these triangles are not preserved between time steps. In fact, they need not even be stored. Unlike the compuation of a triangulation on an entire surface~\cite{Chen1997}, the computation of these local triangles does not require any mapping to a parameter space.

\begin{figure}
  \centering
  \includegraphics[width=0.49\textwidth]{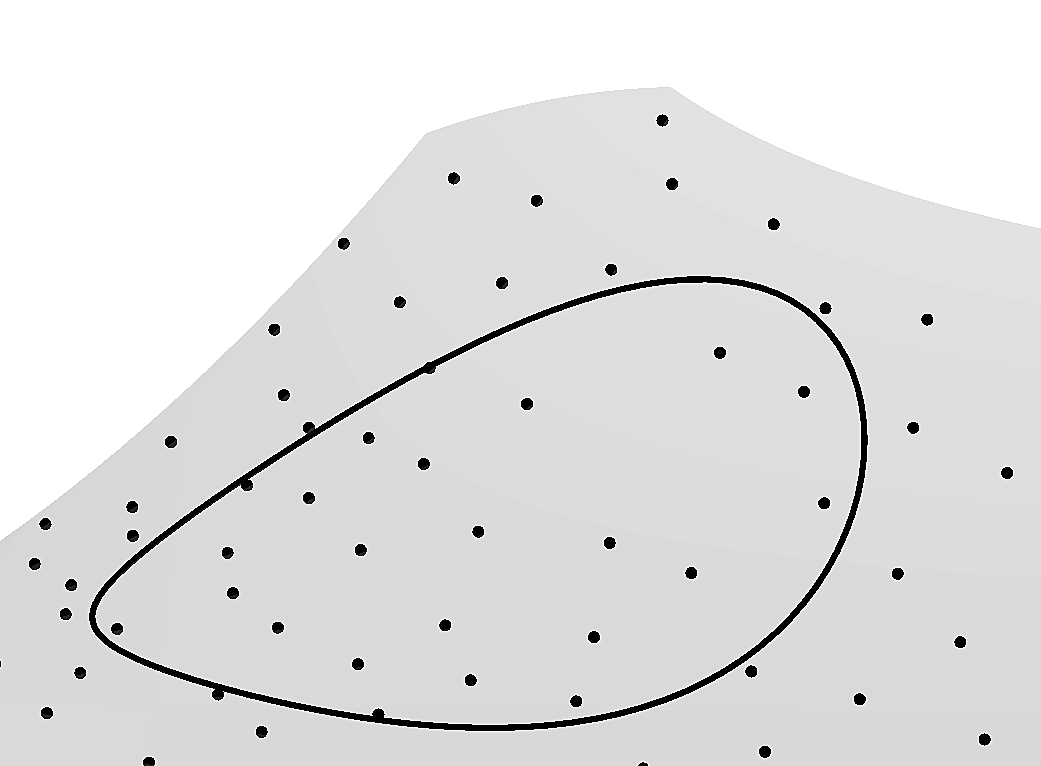}
  \includegraphics[width=0.49\textwidth]{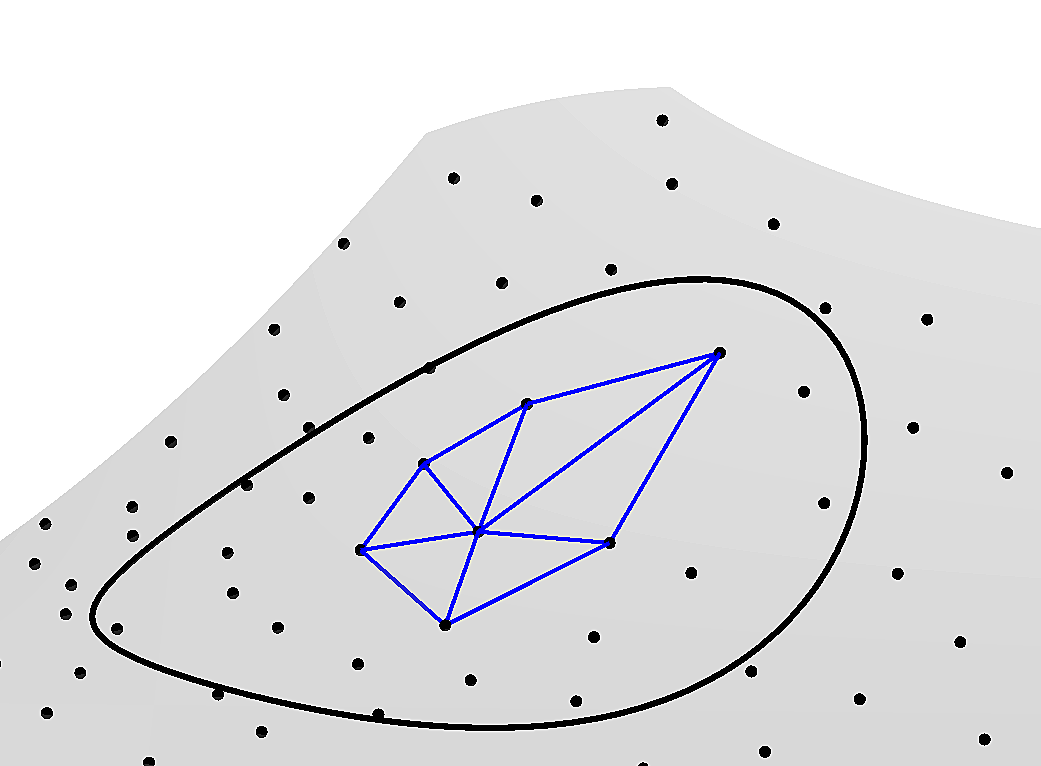}\\
  \includegraphics[width=0.49\textwidth]{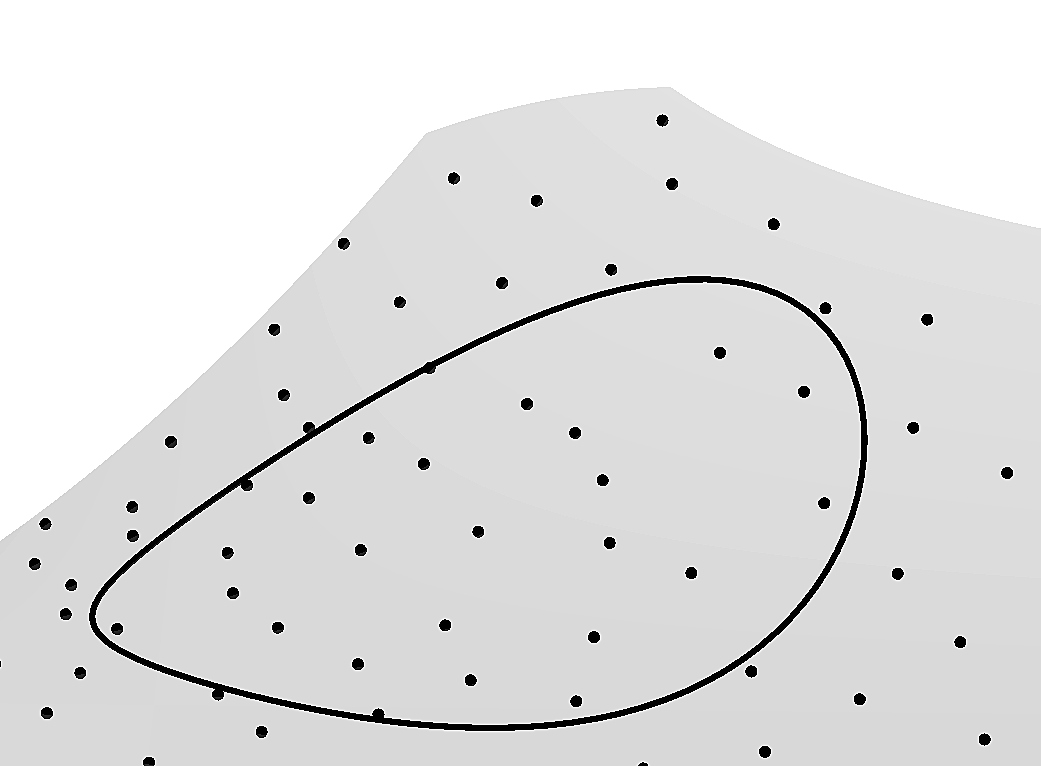}    
  \includegraphics[width=0.49\textwidth]{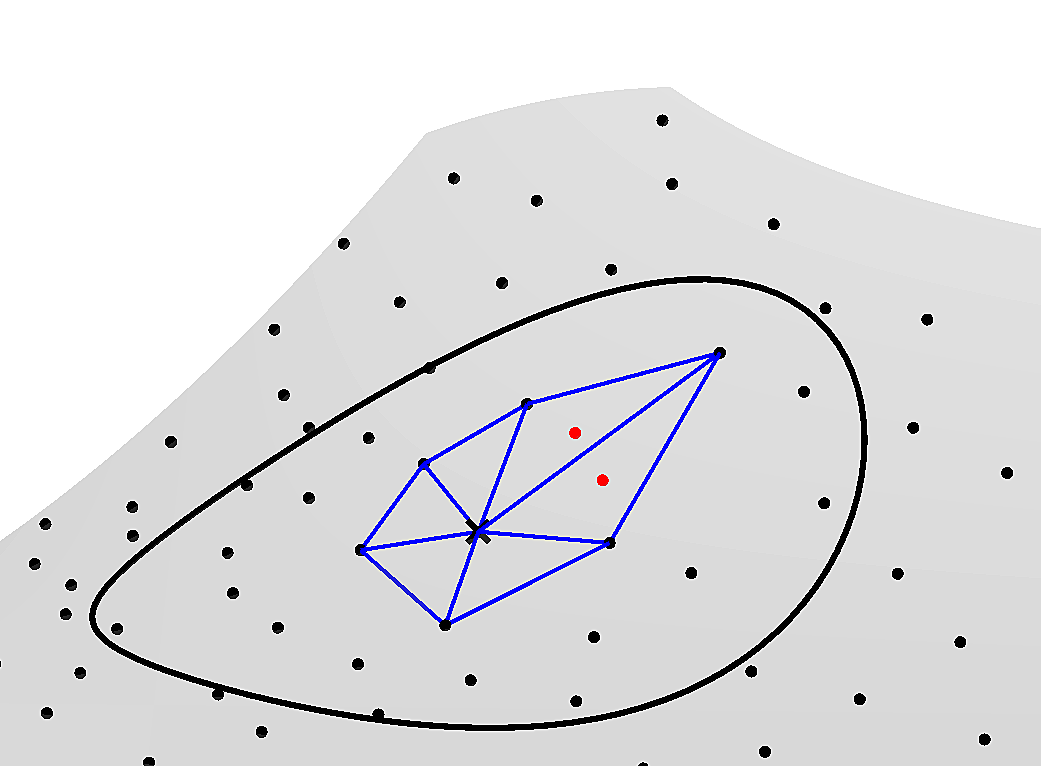}
  \caption{Addition of points in regions with insufficient points. Clockwise from top left: Initial configuration~(top left), triangulation of a point within the support radius~(top right), addition of red points in `large' triangles~(bottom right), and final configuration~(bottom left).}
  \label{Fig:Addition}%
\end{figure}

These locally defined triangles, as shown in Figure~\ref{Fig:LocalTriangulation}, are used to determine locations of holes where points need to be filled. For this, we recall the volumetric point cloud spacing conventions introduced in Section~\ref{sec:Discretization}, which are carried over to the surface case here. We intend to ensure that there is no hole of size $r_{max}h$ in the point cloud, where the smoothing length $h=h(\vec{x}, t)$ can be a function of both space and time. We set $r_{max}=0.45$, following volumetric point cloud spacing conventions~\cite{Drumm2008}. Thus, for any triangle with a circumradius greater than $r_{max}h$, a point is added at its circumcenter, but only as long as the circumcenter is within the triangle. Once this process is carried out at all points in the domain, it is then repeated for the newly created points. This process of filling points is illustrated in Figure~\ref{Fig:Addition}.

For each newly created point, we need to approximate all physical quantities there. A considerable amount of work has been done in determining optimal ways to perform approximations at these new points for the volumetric Lagrangian framework \cite{Drumm2008,Farjoun2008,Iske2007} and those can easily be adapted for the needs of the surface-based case.  

We note that for boundary points, we need to ensure to not form a triangle between three boundary points, as that could lead to addition of points outside the domain. To ensure a sufficient number of points on the boundaries, a separate addition process is also done only for boundary points, with the distance between two adjacent boundary points checked.

%


\subsubsection{Curvature Corrected Addition}~\\
\label{sec:CurvatureCorrectedAddition}
The above process of addition at the circumcenter of large triangles assumes the surface to be piecewise linear. This can be improved by correcting the new location to take the curvature into account. For a point being added, we start by approximating a normal at that point by
\begin{equation}
	\vec{n}_{approx} = \frac{ \vec{n}_1 + \vec{n}_2 + \vec{n}_3}{3}\,,
\end{equation}
where the summation is over the three points in the triangle being considered. The newly created point is then moved in the direction of this approximated normal
\begin{equation}
	\vec{x}_{new} = \vec{x}_c + d_{\kappa} \vec{n}_{approx} \,,
\end{equation}
where $\vec{x}_c$ is the original location at the circumcenter, and $\vec{x}_{new}$ is the corrected location. The distance $d_{\kappa}$ of moving the newly created point is computed based on the angles between the distance vectors and the normal field, in a manner similar to the procedure of normal computation explained earlier. The following constraint is enforced:
\begin{equation}
	\sum_j \delta\vec{x}_j\cdot \vec{n}_j = 0 \,,
\end{equation}
where the sum is over the three points of the triangle, and $\delta \vec{x}_j = \vec{x}_j - \vec{x}_{new}$. This leads to
\begin{equation}
	d_{\kappa} = \frac{ \sum_j ( \vec{x}_j - \vec{x}_c ) \cdot \vec{n}_j}{ \sum_j \vec{n}_{approx}\cdot\vec{n}_j}\,.
\end{equation}

\subsection{Deletion or Merging of Points}
\label{sec:PCMdeletion}

\begin{figure}
  \centering
  \includegraphics[width=0.49\textwidth]{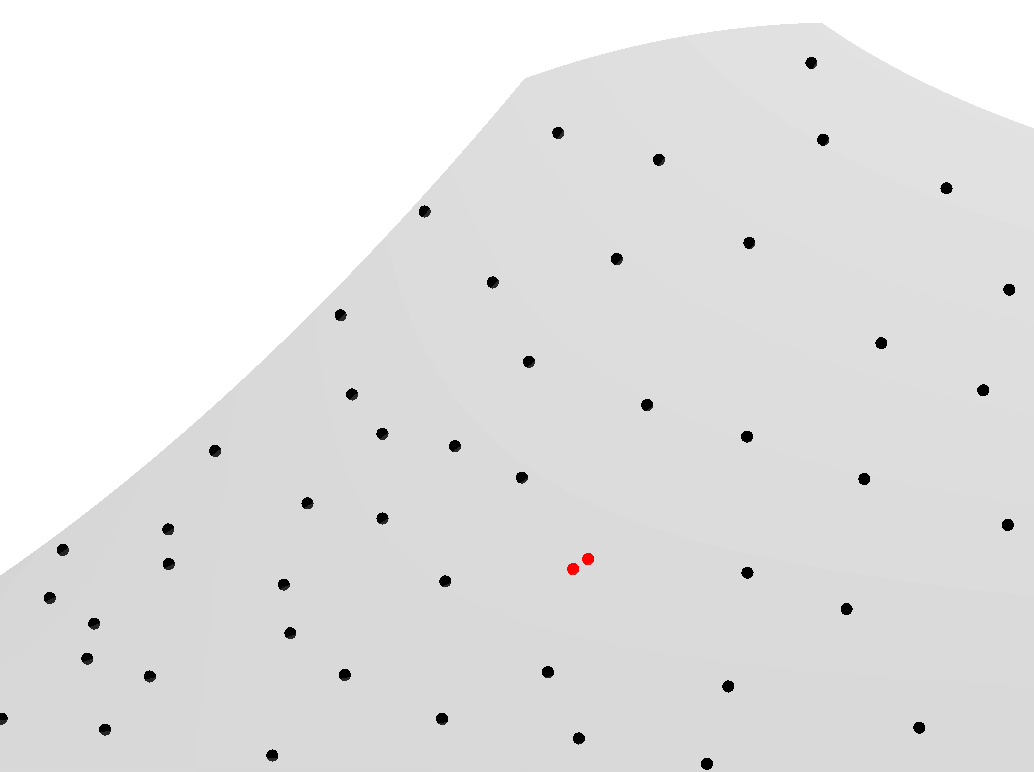}
  \includegraphics[width=0.49\textwidth]{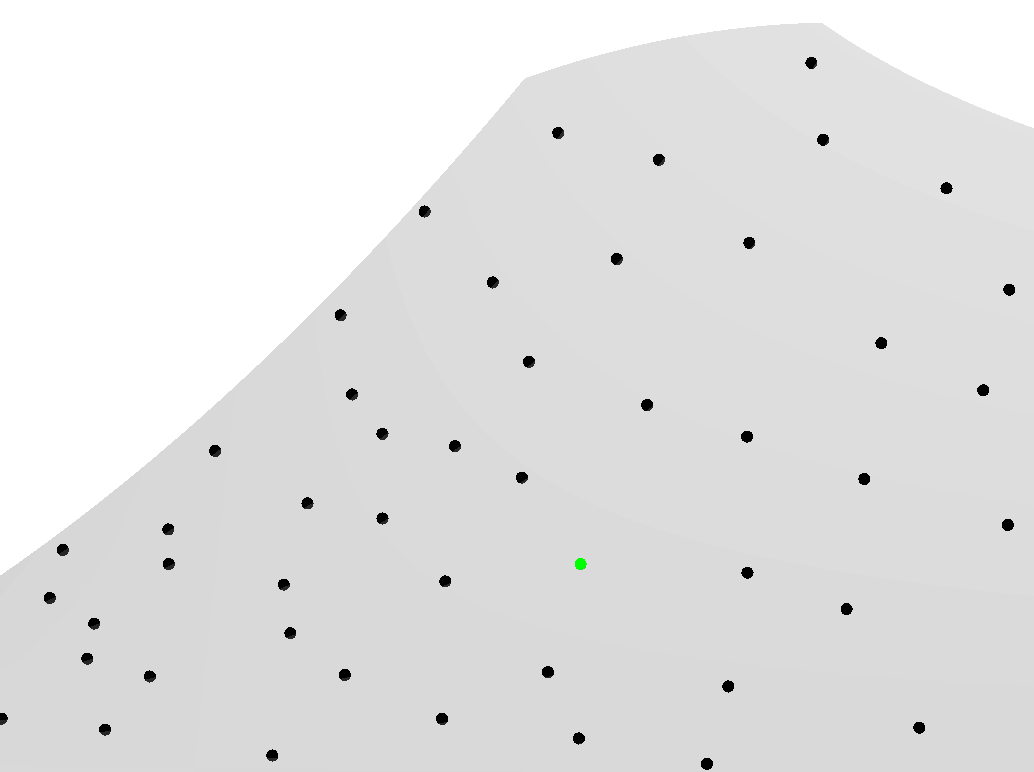}
  \caption{Merging points that are too close. The two red points in the left figure are merged to the green point in the right figure.}
  \label{Fig:Deletion}%
\end{figure}

To prevent clustering of points as a result of Lagrangian movement, numerical points within $r_{min}h$ of each other are merged into a single point, with all physical properties interpolated at the new location. This process is illustrated in Figure~\ref{Fig:Deletion}.

In this process, priority is given to boundary points. If one boundary point and one interior point are within a distance of $r_{min}h$ from each other, instead of merging the two points, the interior point is deleted, and values are re-interpolated at the boundary point. This is done so as to not artificially deform the boundary.

These procedures of addition and deletion of points can be used to obtain adaptive refinement, when needed. Changing the smoothing length $h$ will trigger the required point addition or merging algorithm. For example, refinement can be done based on some error bounds, as has been done in the context of volumetric meshfree GFDMs by \cite{Oanh2017, Urena2017}. Alternatively, refinement can be carried out based on the gradient of a numerical solution $\nabla\phi$. The smoothing length $h$ can be reduced with increasing $\|\nabla\phi\|$, which would result in a higher number of points in those locations. This process is illustrated in a numerical example in Section~\ref{sec:TopoChange}. 

\subsection{Contact Detection and Topological Change}
\label{sec:Contact}

Effectively modelling topological change of manifolds, and contact between different manifolds is one of the most challenging parts of setting up a comprehensive Lagrangian framework. Collision detection is usually done with the help of a mesh \cite{Wriggers2006}. Even particle methods rely on a background mesh for the same \cite{Leung2011}. A lot of work has been done for collision detection between two volume phases, by both mesh-based methods \cite{Wriggers2006} and meshfree methods \cite{Li2001}. Several adaptations of these are needed to apply such methods for collision detection between surfaces given by point clouds.

As a preliminary, each meshfree node is assigned a `chamber' attribute to indicate the difference between different manifolds. Note that for each point $i$, the $n_{S_i}$ points in its neighbourhood $S_i$ all belong to the same chamber. We start with collision between different manifolds, and later extend these algorithms for self-intersection which involves contact between different parts of the same manifold. 

In this paper, we do contact detection in two steps. For a particular point $i$ in chamber $c_i$, the first step involves a spatial search to determine which points of other chambers could possibly come into contact with the current chamber. To do this efficiently, and to avoid the need to compute distances with every point of other chambers, we use the neighbour searching algorithms mentioned above. After this stage, we assume that each point $i$ has a list of own-chamber neighbours $S_i$, and a list of other-chamber `neighbours' $S_i^{ext}$ which are within a distance $h$ of point $i$. Note that, unlike $S_i$, there is no need to store $S_i^{ext}$ after contact detection is done. Furthermore, $i\in S_i$ but $i \notin S_i^{ext}$.

The next step involves the actual contact determination. We propose doing this in two parts. One involves determining penetrated points, i.e. determining parts of a surface that have `crossed' another surface. For this, we first compute a signed distance $DC_i$ to the other chamber, where $DC_i$ indicates the distance to contact. $DC_i$ is computed as the distance of the point $i$ from an approximated surface given by the points in $S_i^{ext}$. The sign depends on the normals of $S_i^{ext}$. First, the distance to the surface given by $S_i^{ext}$ is approximated as a quadratic polynomial
\begin{equation}
	\label{Eq:DTB}
	d^{S_i^{ext}} = d_0 + d_1x + d_2y + d_3x^2 + d_4y^2 + d_5xy\,,
\end{equation}
where the coefficients $d_k$ are determined based on the positions and normals of the points in $S_i^{ext}$. For this, we use the fact that the distance should be $0$ at every point in $S_i^{ext}$. Thus,  $d^{S_i^{ext}}( {\vec{x}_m} ) = 0\,,\, \forall m \in S_i^{ext}$. Furthermore, the gradient of the distance function should be along the normal to the surface. At a dummy location obtained by moving a point in $S_i^{ext}$ by a fixed distance $\xi$ in the normal direction, the distance function should return $\xi$. Thus, $d^{S_i^{ext}}( {\vec{x}_m} \pm \xi \vec{n}_m ) = \pm \xi \,,\, \forall m \in S_i^{ext}$ and sufficiently small $\xi>0$. We take $\xi$ to be a third of the minimum distance between points. These two conditions together form, in general, an over-determined system, which is solved in a least squares sense to determine the coefficients $d_k$. Then, $DS_i$ is computed by evaluating Eq.\,\eqref{Eq:DTB} at $\vec{x}_i$. 

\begin{figure}
  \centering
  \includegraphics[width=0.4\textwidth]{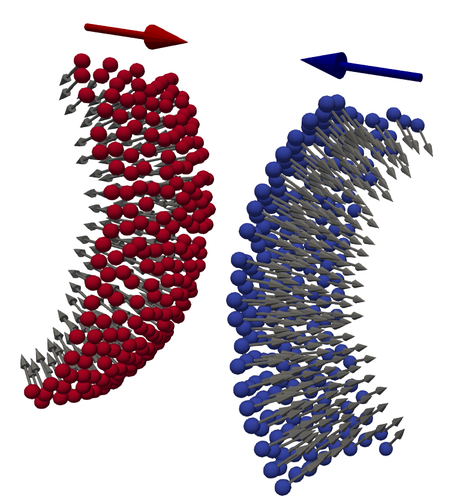}
  \hspace{50pt}
  \includegraphics[width=0.4\textwidth]{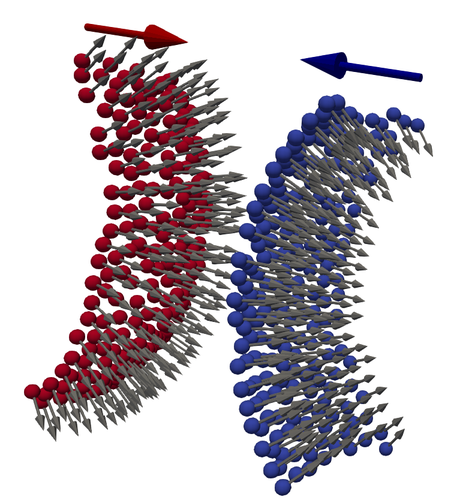}
  \caption{Different possible orientations of manifolds coming into contact. The colours of the points indicate the different manifolds. The grey arrows are the manifold normals. The red and blue arrows indicate the velocities of the manifold. The blue and red manifolds are moving towards each other in both cases. In the figure on the left, the normals of the two manifolds point away from each other. While in the figure on the right, they point towards each other.}
  \label{Fig:NormalDirection}%
\end{figure}

It is important to note that, unlike the volume-based case, the sign of the distance computed above does not directly determine whether or not the point has penetrated another chamber. In the volume case, each boundary surface has normals pointing either inwards or outwards, which is known a-priori. However, in the case of manifolds, the concepts of interior and exterior are not always well defined. And thus, the direction of the normals are not always known a-priori. This is illustrated in Figure~\ref{Fig:NormalDirection}. In Figure~\ref{Fig:NormalDirection}, the red and blue manifolds are moving towards each other. It shows the two possibilities of the normal direction of the red manifold. As shown in the figure on the right, it is quite possible that the normals of the chambers are pointing in the same direction. Thus, checking for conflicting Lagrangian information, as is often done \cite{Petras2016}, is not always enough. 

Thus, the sign of $DC_i$ at one time step alone does not give us enough information. To solve this problem, we store the information about signed distance for each point. If the sign of $DC_i$ has flipped compared to the previous time step, then we can say that the point has penetrated another chamber. Using a soft penetration model, $|DC_i|$ is then used for force computation, which depends on the considered physical model, as would be done in mesh-based methods \cite{Wriggers2006}.

The above method checks for opposing chambers that crossed each other in the previous movement step. The second part of contact detection involves determining opposite chamber points that are almost in contact. This is done by keeping track of points that are within the deletion distance of $r_{min}h$, but have not yet crossed into the opposite chamber. Once again, force computation can be done based on the distance, which will depend on the considered physical model.

The processing of detected contacts between manifolds is differentiated into $3$ types:

\begin{itemize}
\item \textit{Non-penetration contact:}  This is the most standard form of contact, indicating that two surfaces can not cross each other. Handling this numerically takes the form of adding the correct force once particles penetrate or come too close to another chamber, as described above.

\item \textit{Delete contact:} This is the case, for example, when two fluid droplets modelled only by their surfaces come into contact with one another. Thus, the points of each chamber coming into contact are deleted. For the geometric handling of this case, a numerical point is deleted if it has just penetrated into another chamber, or if it is within deletion distance of a point of another chamber.

\item \textit{Merge contact:} This third situation happens when two manifolds `stick' together, such that the resultant can be treated as a single manifold, with modified properties. For this, points of different chambers are merged into a single chamber. 

\end{itemize}

Numerical examples of the first two types of contact, based purely on the geometry of the problem, without any physical model, are presented in Section~\ref{sec:TopoChange}. 

The above methods can also be extended quite easily for contact detection between different parts of the same manifold. For this, the only difference is the identification of the equivalents for $S_i$ and $S_i^{ext}$ defined above, which is done based on the direction of the normals.

\subsection{Numerical Results}
\label{sec:NumRes}

Before solving PDEs on evolving surfaces, we first test the Lagrangian framework by only performing advection of the manifold. In each of the simulations presented below, the point clouds are irregularly spaced. The initial point clouds are obtained starting from a CAD file by placing points using an advancing front technique for point clouds, similar to that done by \cite{Drumm2008}.

\subsubsection{Rotating Quarter Sphere}\label{sec:RotatingQSphere}
\begin{figure}
  \centering
  \includegraphics[width=0.6\textwidth]{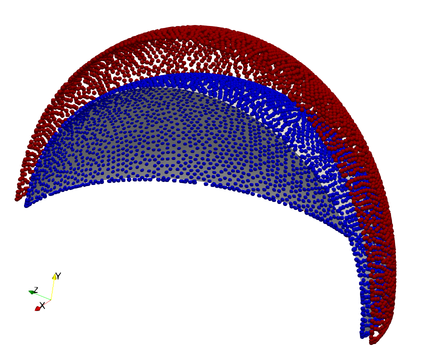}
  \caption{Quarter sphere point clouds after one full rotation. The blue point cloud is a result of the second order movement, Eq.\,\eqref{Eq:Move_o2}, while the red point cloud is a result of the first order movement, Eq.\,\eqref{Eq:Move_o1}. The grey background shell represents the original configuration at initial time.}
  \label{Fig:O1vsO2}%
\end{figure}

~\\
To show the need of movement by the second order method according to Eq.\,\eqref{Eq:Move_o2}, we consider the case of a quarter of a unit sphere rotating about its center. To illustrate the effect of the movement only, no addition and deletion of points is done in this example. The sphere is centered at the origin, and the velocity of the surface is given by
\begin{equation}
	\vec{v} = \left( \begin{array}{c}
		y\\
		-x\\
		0
	\end{array} \right) + \left( \begin{array}{c}
		0\\
		z\\
		-y
	\end{array} \right)\,,
\end{equation}
where the terms correspond to a rotation about the $z$ and $x$ axis, respectively. Clearly, the velocity is tangential to the surface at every point. Thus, updating point locations with the conventionally used first order movement, Eq.\,\eqref{Eq:Move_o1}, results in points going off the sphere. As a result, numerically, the radius of the quarter sphere increases with every time step. This issue is no longer present for the second order movement, Eq.\,\eqref{Eq:Move_o2}. Figure~\ref{Fig:O1vsO2} shows the point clouds after one full rotation based on the two different movement methods in comparison with the original quarter sphere. The quarter sphere surface is discretized with $N=3\,792$ irregularly spaced points. The time step is fixed as $\Delta t = 0.05$. Figure~\ref{Fig:O1vsO2} illustrates the inaccuracy of the first order movement whenever the advection velocity has a rotational component. The blue point cloud from the second order movement is virtually at the initial location, while the red one from the first order movement is off it. To quantify this error, we measure the distance from the unit sphere. After one full rotation, the mean distance of the resultant point clouds from the unit sphere is $1.37\times 10^{-1}$ for the first order movement, and $5.57\times 10^{-4}$ for the second order movement. It is important to note here that for the first time step, the second order method reverts to the first order one due to unavailability of a previous velocity, and that the mean error after the first time step is $1.12\times 10^{-4}$ for both methods. The error as a result of the second order movement barely increases after the first time step, whereas that of the first order movement increases by several orders of magnitude. Thus, we now only use the movement according to Eq.\,\eqref{Eq:Move_o2} for the remainder of the paper.

\subsubsection{Deforming Hemisphere}
\begin{figure}
  \centering
  \includegraphics[width=0.32\textwidth]{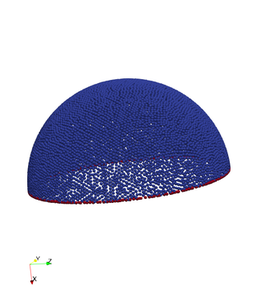}
  \includegraphics[width=0.32\textwidth]{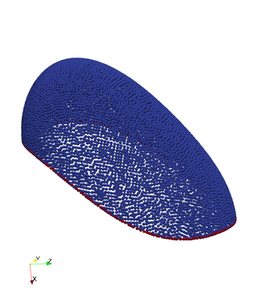}
  \includegraphics[width=0.32\textwidth]{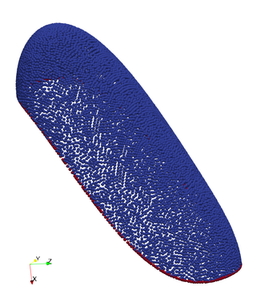}
  %
  \includegraphics[width=0.32\textwidth]{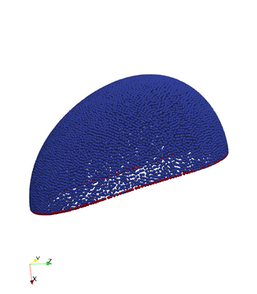}
  \includegraphics[width=0.32\textwidth]{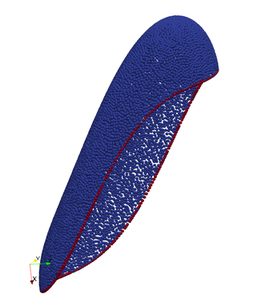}  
  \includegraphics[width=0.32\textwidth]{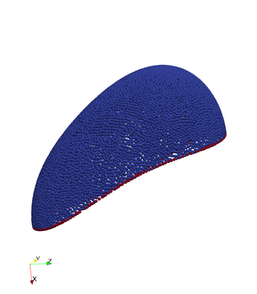}
  \caption{Elongation and contraction of a unit hemisphere, with addition and deletion of points to maintain regularity. Clockwise from top left: At times $t = 0$~(top left), $t = 0.09$~(top center), $t = 0.24$~(top right), $t = 0.55$~(bottom right), $t = 0.745$~(bottom center) and $t = 0.975$~(bottom left). The red points indicate manifold boundaries.}
  \label{Fig:HemisphereStretch}%
\end{figure}

\begin{figure}
  \centering
  \includegraphics[width=0.8\textwidth]{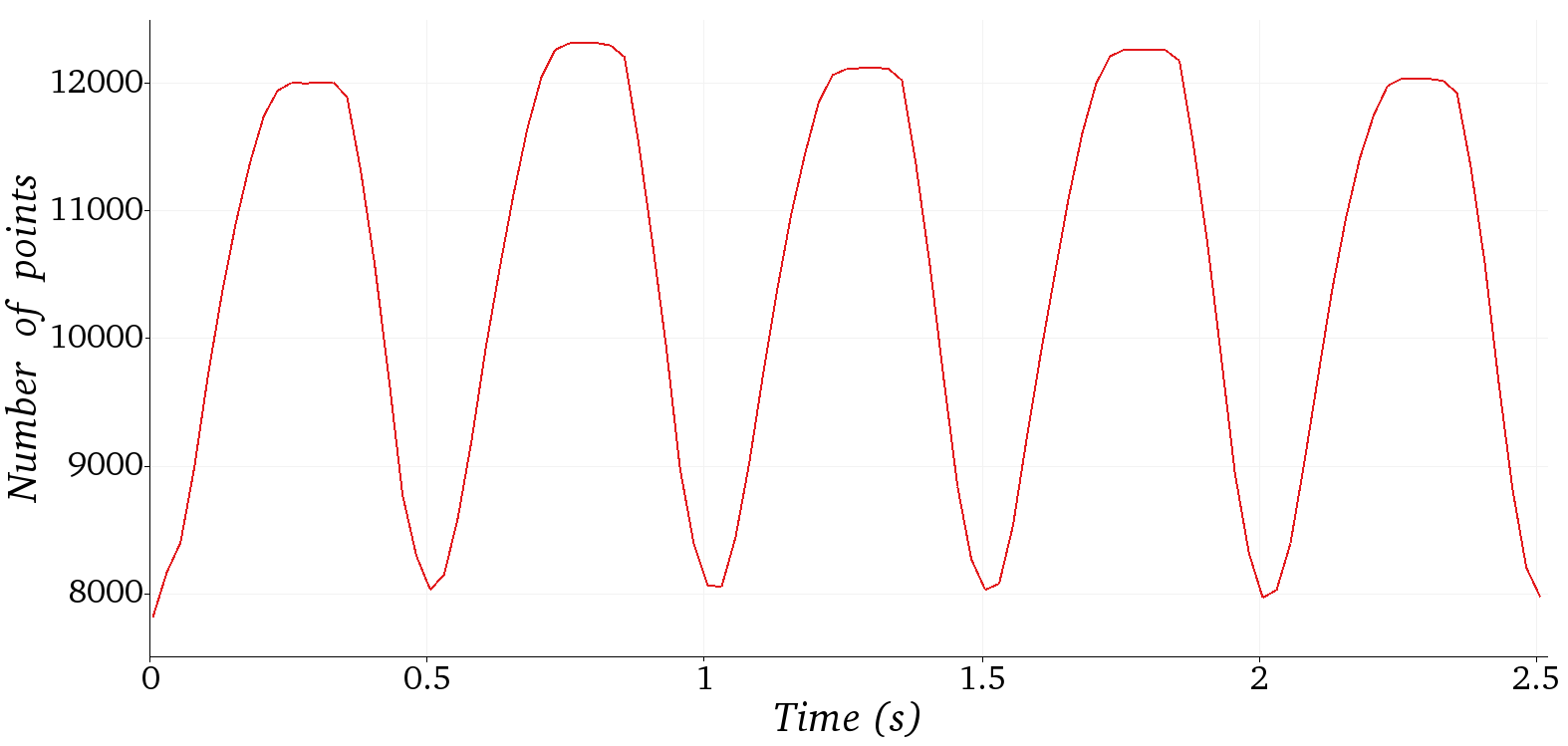}
  \caption{Number of points as a function of time for the distorting hemisphere test case.}
  \label{Fig:NumberOfPoints}%
\end{figure}

~\\
To illustrate the addition and deletion of points in the Lagrangian framework, we consider stretching and contracting of a unit hemisphere as shown in Figure~\ref{Fig:HemisphereStretch}. The velocity field is given by
\begin{equation}
\label{Eq:HemisphereStretchVel}
 \vec{v} = \left( \begin{array}{c}
 			2\pi \cos(2\pi t) \sin(\frac{\pi}{2}z) \\
 			0\\
 			0
	 	 \end{array} \right) \,,
\end{equation}
which is a modification of the example considered in \cite{Chen2016}. 
Figure~\ref{Fig:HemisphereStretch} shows the evolution of the surface for $\Delta t = 0.005$. The number of points $N$ on this evolving surface as a function of time is shown in Figure~\ref{Fig:NumberOfPoints}. As the sphere is stretched, to maintain regularity of the point cloud, points are added on the manifold. For the same reason, points are deleted from the manifold as the sphere contracts. 

\begin{figure}
  \centering
  \includegraphics[width=0.45\textwidth]{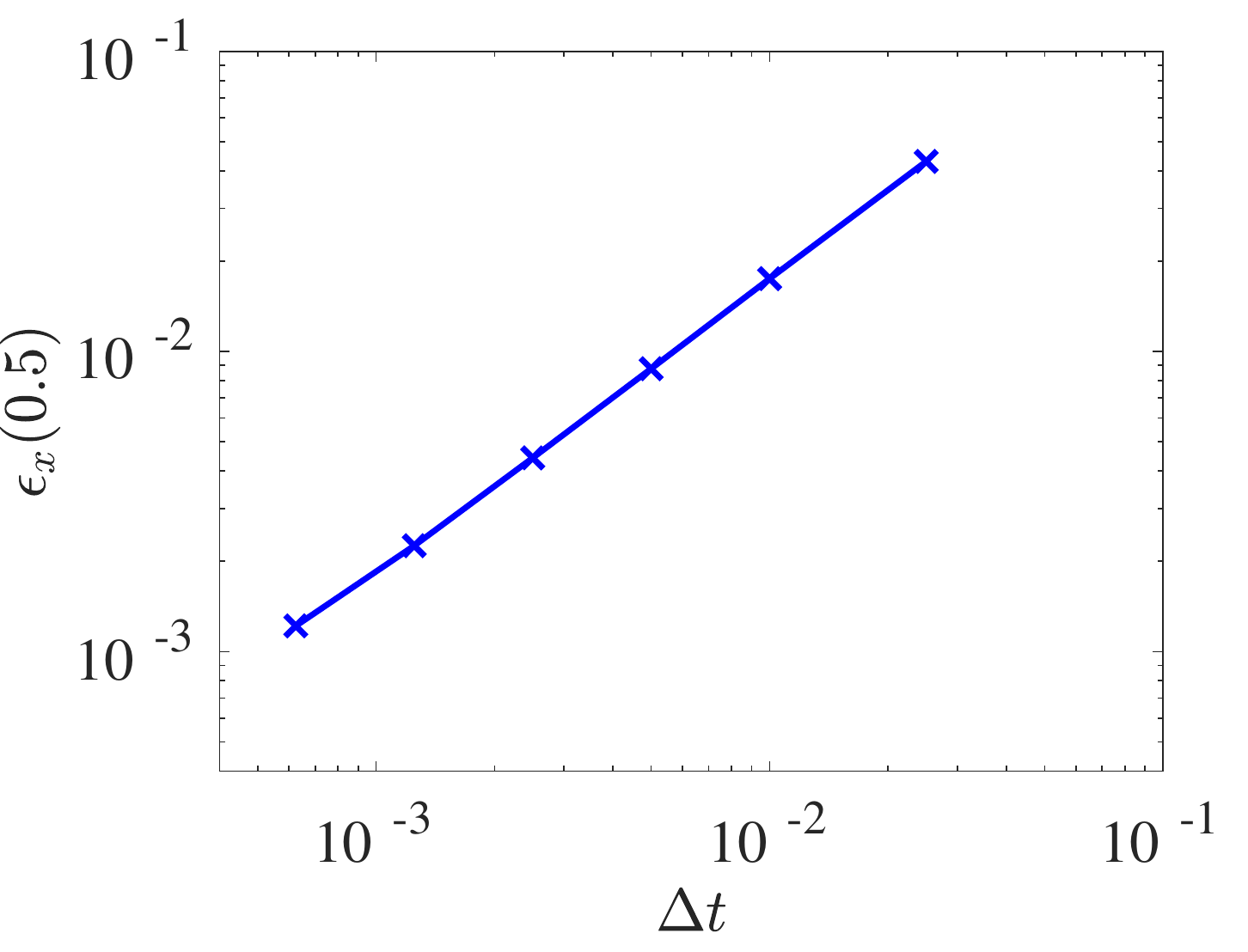}
  \includegraphics[width=0.45\textwidth]{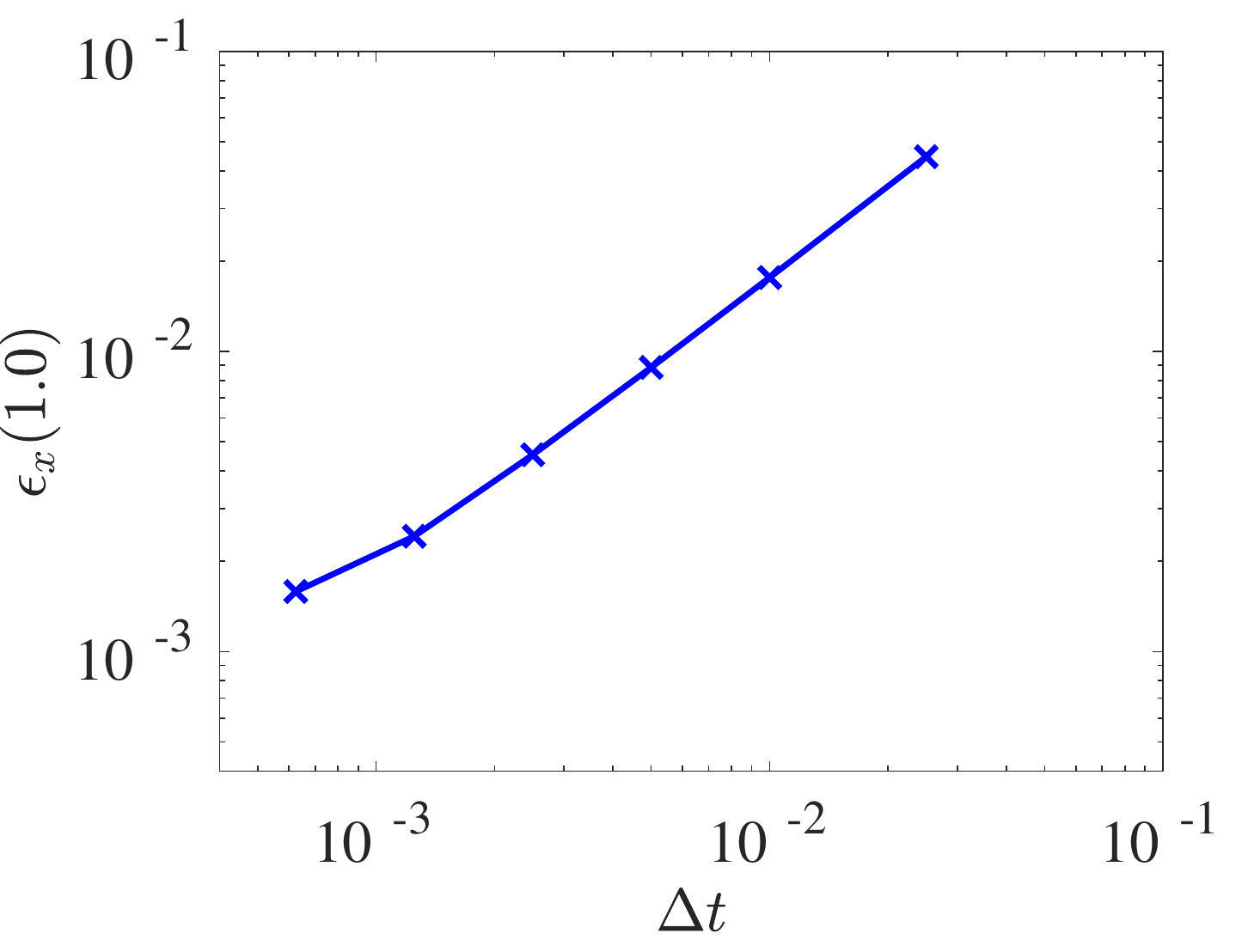}
  \caption{Convergence of error in location as the deforming hemisphere returns to its original state at $t=0.5$~(left) and $t=1.0$~(right).}
  \label{Fig:HemisphereStretch_Error}%
\end{figure}

This serves as a good example of surface deformation because the surface returns to its original shape. We know that at every multiple of $t=0.5$, the manifold should once again take the initial shape. The resulting error in the Lagrangian movement can thus be checked by measuring the mean distance from the unit sphere, as done in the previous example.
\begin{equation}
	\label{Eq:ElongatingHemisphereError}
	\epsilon_{x}(t) = \frac{1}{N(t)} \sum_{i=1}^{N(t)}  \left| \| \vec{x}_i \|^2 - 1  \right| \,.
\end{equation}
We reiterate that $\epsilon_x$ should be $0$ only for integer multiples of $t=0.5$. The convergence of the error $\epsilon_x$ at $t=0.5$ and $t=1.0$ is plotted in Figure~\ref{Fig:HemisphereStretch_Error} for a varying time step. The two plots also show that the accumulation of the Lagrangian movement error with time is not significant. The figure illustrates that a numerical convergence order of $\mathcal{O}(\Delta t)$ is observed. The main reason for the small convergence order is that the rate determining step is the very first time step, where a first order movement is necessary as no previous velocity is available. We note that the same was quantified in the previous example. This could be avoided by prescribing an initial velocity derivative, which is not done here.

%
%
%

\subsubsection{Joining Spheres}\label{sec:TopoChange}~\\
To illustrate the contact handling algorithms, we consider two unit spheres moving towards each other. Simulations are done for the same with both delete contact, and non-penetration contact enforced, as described in Section~\ref{sec:Contact}.

The two examples considered in this section are only used to illustrate the applicability of the geometric collision detection and penetration avoidance algorithms presented earlier. Thus, no physical model has been applied. The delete contact case is handled as described in  Section~\ref{sec:Contact}. For the non-penetration contact case, since no force is added, a slight modification is required to avoid penetration. Once penetration has been detected geometrically, as described earlier, penetrated points are projected back so as to ensure that the two discrete manifolds do not cross each other. Velocities are recomputed as averaged values after the projection.

In both simulations, we use a time step of $\Delta t = 0.03$, and a smoothing length of $h=0.1$ which corresponds to $N = 15\,018$ total points in each of the two spheres at the initial time level. The spheres move towards each other with a velocity of $\vec{v}=(-\text{sign}(x_1)  0.5, 0, 0 )$, where `sign' is the signum function, and the origin is between the two spheres at initial state. The center of the spheres are at an initial distance of $2.2$ units apart. 

\begin{figure}
  \centering
  \includegraphics[width=0.49\textwidth]{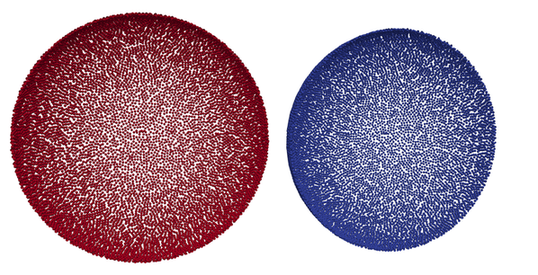}
  \includegraphics[width=0.49\textwidth]{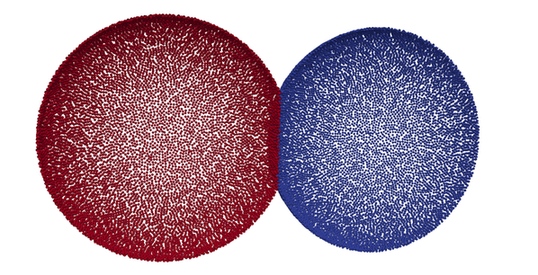}
  %
  \includegraphics[width=0.49\textwidth]{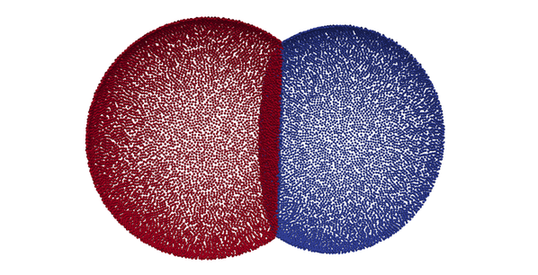}  
  \includegraphics[width=0.49\textwidth]{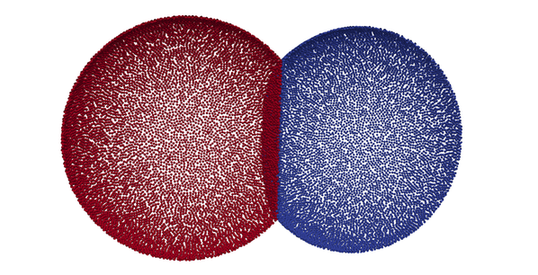}
  %
  \includegraphics[width=0.49\textwidth]{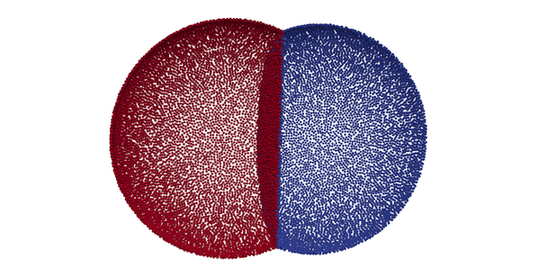}  
  \includegraphics[width=0.49\textwidth]{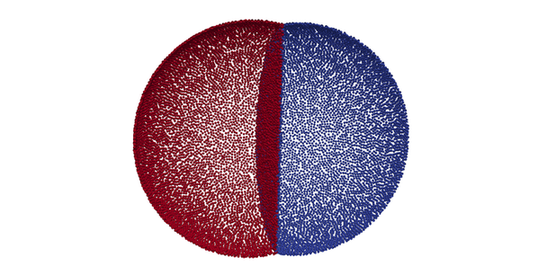}
  \caption{Two sphere surfaces moving towards each other, with non-penetration contact enforced. Clockwise from top left: At times $t = 0$~(top left), $t = 0.36$~(top right), $t = 0.72$~(middle right), $t = 1.08$~(middle left), $t = 1.44$~(bottom left) and $t = 1.8$~(bottom right). The colour represents the different manifolds. To make the result easier to visualize, only the half of the domain away from the viewing angle is shown in the figure.}
  \label{Fig:NonPenetrationContact}%
\end{figure}

For the non-penetration contact case, Figure~\ref{Fig:NonPenetrationContact} shows the evolution of the spheres as they move towards one another. Note that Figure~\ref{Fig:NonPenetrationContact} only shows a slice of the two spheres, such that only the half away from the viewing angle is seen.

\begin{figure}
  \centering
  \includegraphics[width=0.49\textwidth]{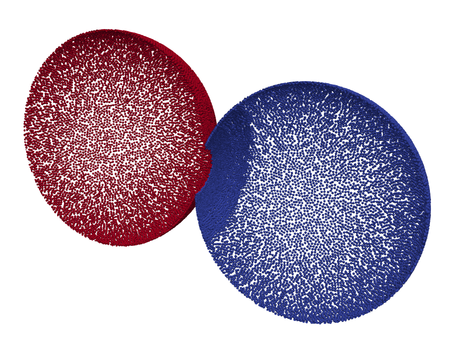}
  \includegraphics[width=0.49\textwidth]{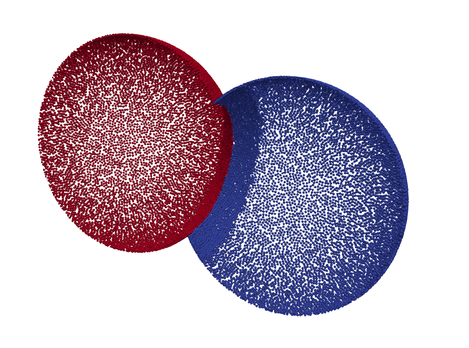}
  %
  \includegraphics[width=0.49\textwidth]{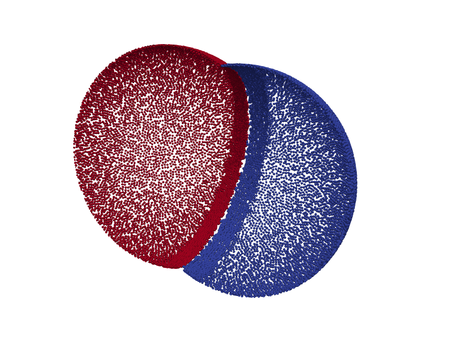}  
  \includegraphics[width=0.49\textwidth]{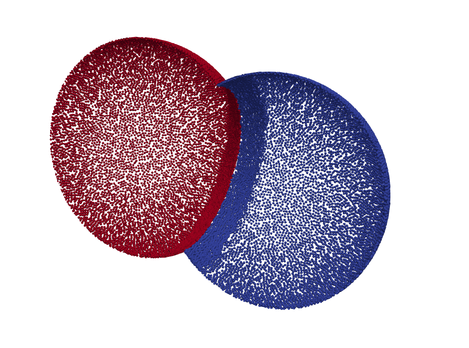}  
  \caption{Two sphere surfaces moving towards each other, with delete contact enforced. Clockwise from top left: At times $t = 0.03$~(top left), $t = 0.63$~(top right), $t = 0.93$~(bottom right), and $t = 1.5$~(bottom left). The colour represents the different manifolds. To make the result easier to visualize, only the half of the domain away from the viewing angle is shown in the figure.}
  \label{Fig:DeleteContact}%
\end{figure}

The delete contact simulation uses the same setup with one difference. To better process the contacts, we apply an adaptive refinement. For regions closer to the interface, the smoothing length $h$ is reduced to $h=0.05$. As a result, the hole filling algorithm explained in Section~\ref{sec:PCMaddition} causes the point cloud to become finer near the interface. Figure~\ref{Fig:DeleteContact} shows the result when enforcing the delete contact. It also illustrates the adaptive refinement. Once again, only a slice of the result is shown. We note that the jagged interface between the two surfaces after deletion on contact is a result of unevenly spaced points. This becomes more regular (more straight, in this case) as the number of points is increased, or as refinement is done near the interface.

\subsubsection{Addition of Points}~\\
The examples considered so far use addition of new points at the circumcenter of `large' triangles. We now consider curvature corrected addition as per Section~\ref{sec:CurvatureCorrectedAddition}. We consider the case of addition of points on the stationary surfaces of a sphere, and of a torus given by
\begin{equation}
	\left( c - \sqrt{x^2 + y^2} \right)^2 + z^2 = a^2\,,
\end{equation}
with $c=3$ and $a = 1$. In each case, the smoothing length is decreased which results in the point cloud getting finer. The error in the location of addition of new points is measured as the mean distance off the surface. This error is give by Eq.\,\eqref{Eq:ElongatingHemisphereError} for the sphere. For the torus, it is given by
\begin{equation}
	\epsilon_{x} =  \frac{1}{N(t)} \sum_{i=1}^{N(t)}  \left| \left( c - \sqrt{x^2 + y^2} \right)^2 + z^2 - a^2 \right| \,.
\end{equation}

\begin{figure}
  \centering
  \includegraphics[width=0.45\textwidth]{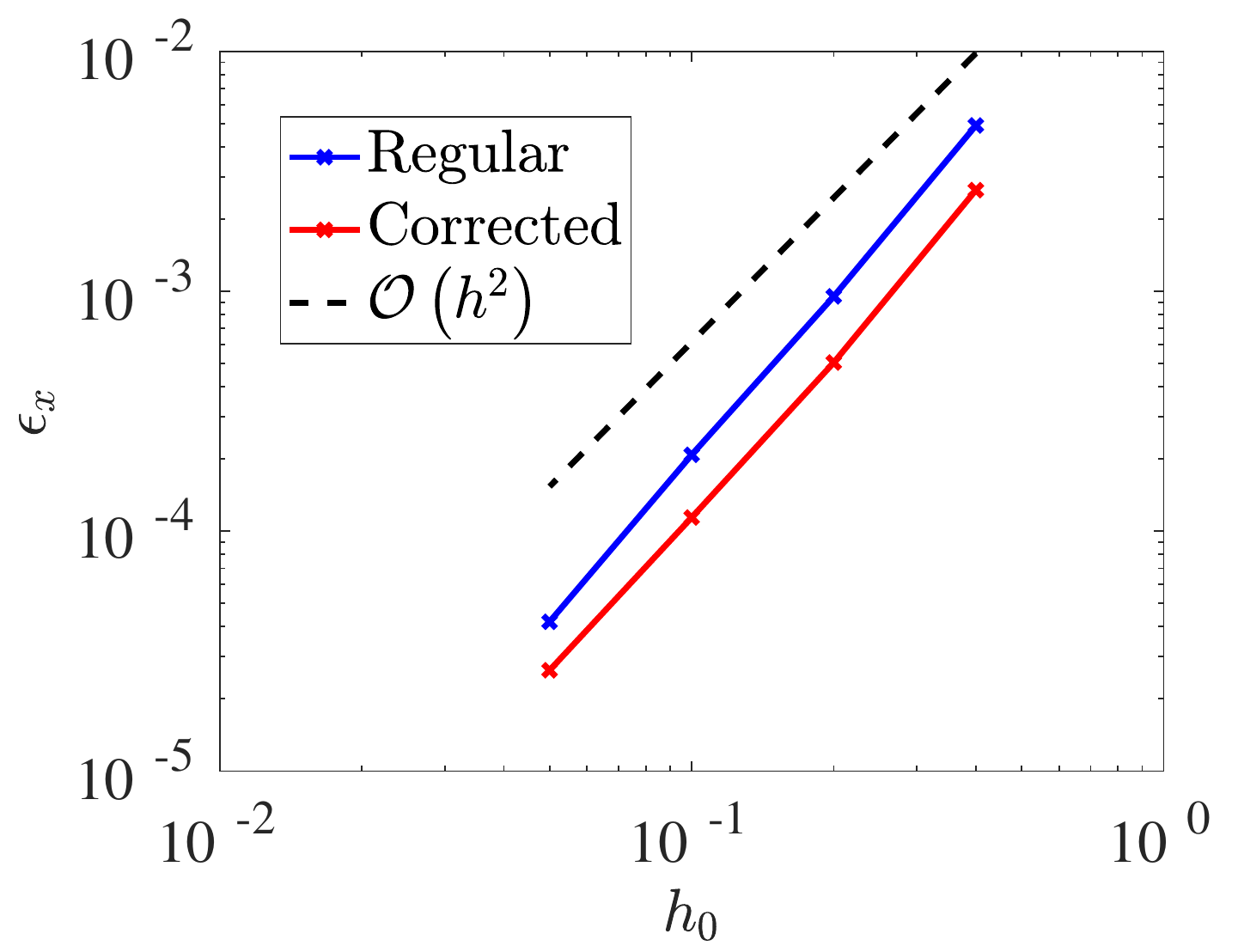}
  \includegraphics[width=0.45\textwidth]{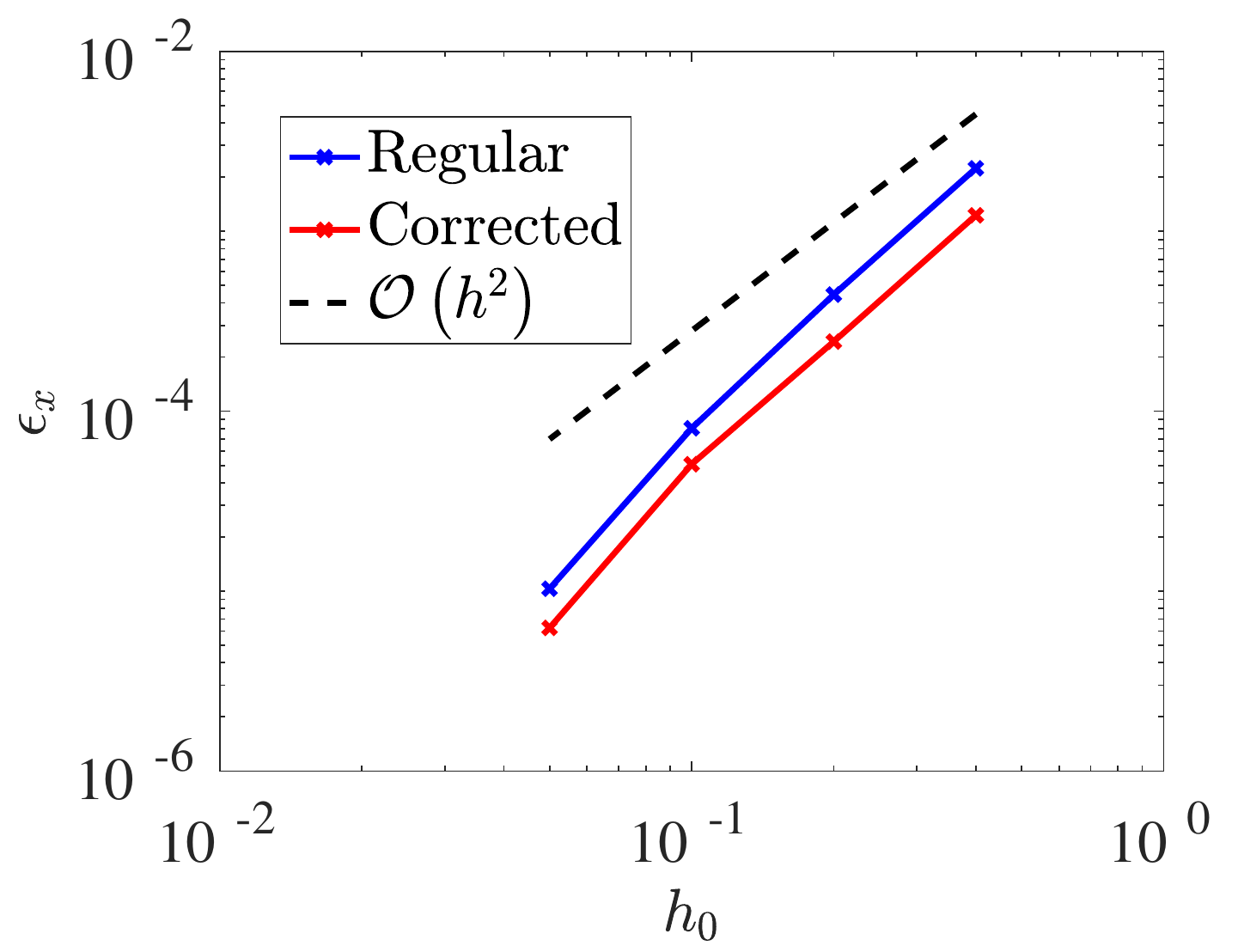}
  \caption{Error in location of newly added points against the discretization size. The blue lines indicate the errors after regular addition at the circumcenter. The red lines indicate the errors when the curvature corrected addition according to Section~\ref{sec:CurvatureCorrectedAddition} is used. Errors are shown for the surface of a sphere~(left) and the surface of a torus~(right).}
  \label{Fig:AdditionError}%
\end{figure}

For both considered geometries, points are added to $4$ point clouds of different discretization sizes to observe how the error in location converges with increasing number of points. In each case, the initial smoothing length $h_0$ is halved, which results in approximately quadrupling the number of points. The resultant error after the entire addition procedure is plotted in Figure~\ref{Fig:AdditionError} for different starting discretizations. The corresponding number of points in each domain before the refining is done is shown in Table~\ref{Tab:AdditionErrors_NumOfPoints}. Figure~\ref{Fig:AdditionError} illustrates that using the curvature corrected addition results in about halving the error, and that the errors converge faster than $\mathcal{O}(h^2)$. Further improvements to curvature corrected addition remains an open problem, and will not be discussed in the present work.

\begin{table}[!htbp]
\caption{Number of points~(before refining) corresponding to the discretization sizes mentioned in Figure~\ref{Fig:AdditionError}. In each case, the final number of points after refining, when the errors are measured, are approximately $4$ times that of the initial number mentioned below.}
\centering
\label{Tab:AdditionErrors_NumOfPoints}
	\begin{tabular}{|l|r|r|}
\hline
	$h_o$ & Sphere & Torus \\
\hline
&& \\[\dimexpr-\normalbaselineskip+2pt]
	$0.4$  & $770$  & $7\,801$\\
	$0.2$  & $3\,520$ & $34\,836$\\
	$0.1$ & $14\,994$ & $147\,112$\\
	$0.05$& $63\,627$ & $612\,828$\\
\hline	
	\end{tabular}
\end{table}
%


\section{PDEs on moving manifolds}
\label{sec:PDEsMovingMan}

We now use the meshfree Lagrangian framework for surfaces developed above to solve PDEs on moving manifolds. In contrast to many existing methods for PDEs on evolving surfaces \cite{Chen2016}, in this method, integration is done on the surface at the existing time level, and no map is maintained between the initial and current state. At each time step, differential operators are computed for the point positions at that level directly. 

The moving Lagrangian framework for evolving PDEs developed in Section~\ref{sec:LagrangianFramework} of this paper can be used to solve PDEs on moving manifolds with a wide variety of methods to define the numerical derivatives. In this paper, we use meshfree Generalized Finite Difference Methods~(GFDMs) for surface PDEs for the same. This is based on our earlier work for PDEs on stationary surfaces~\cite{Suchde2019}. For the sake of completeness, we provide a short description about how they are computed. 

\subsection{Differential Operators}
\label{sec:DiffOp}

Meshfree GFDMs \cite{Fan2018,Gavete2017,Katz2010,Luo2016} are strong form methods. For volume-domains, they have been shown to be robust methods, with a wide variety of applications \cite{Drumm2008,Jefferies2015,Michel2017,Moller2007}. They have also been referred to under the name of Generalized Moving Least Squares~(GMLS) \cite{Mirzaei2012}. Differential operators are computed using a least squares approach while ensuring that monomials up to a certain order are differentiated exactly.

Here, we use the GFDM formulation for surface PDEs from \cite{Suchde2019}. This relies on discretizing differential operators entirely in the tangent space. This method scales with the true dimension of the manifold. It has the advantage of avoiding differential geometry complexities. One of the biggest advantages of this method is that most existing developments in volume-based numerical methods can directly be carried over to surfaces. GFDMs also have the advantage that it is straightforward to handle a wide variety of boundary conditions. 

Differential operators are computed by virtual projections to the tangent space at each point which recovers a Euclidean problem in the dimension of the tangent space. Below, we explain how the differential operators are computed. For details about the same, we refer to \cite{Suchde2019}.

\subsubsection{Surface Gradients}~\\
For a function $u$ defined on the surface, its surface gradient is given by
\begin{equation}
	\nabla_{M} u = \nabla \hat{u} \,,
\end{equation}
where $\nabla_M$ is the surface gradient, $\nabla$ is the volumetric gradient, and $ \hat{u}$ is a normal extension of $u$ to a band around the surface such that $\vec{n}\cdot\nabla \hat{u}= 0$. Thus, computing a numerical approximation to $\nabla_{M} u$ can be done by computing a numerical approximation to $ \nabla \hat{u}$. Since $\vec{n}\cdot\nabla \hat{u}= 0$, we only have to approximate the tangential components $\vec{t}_k \cdot \nabla \hat{u}$ for $k=1,2$. For each point $i$, these tangential derivatives $\vec{t}_k \cdot \nabla \hat{u}$ are approximated in the tangent plane $T_i$ spanned by $\vec{t}_{1,i}$ and $\vec{t}_{2,i}$. 

\begin{figure}
  \centering
  \includegraphics[width=0.5\textwidth]{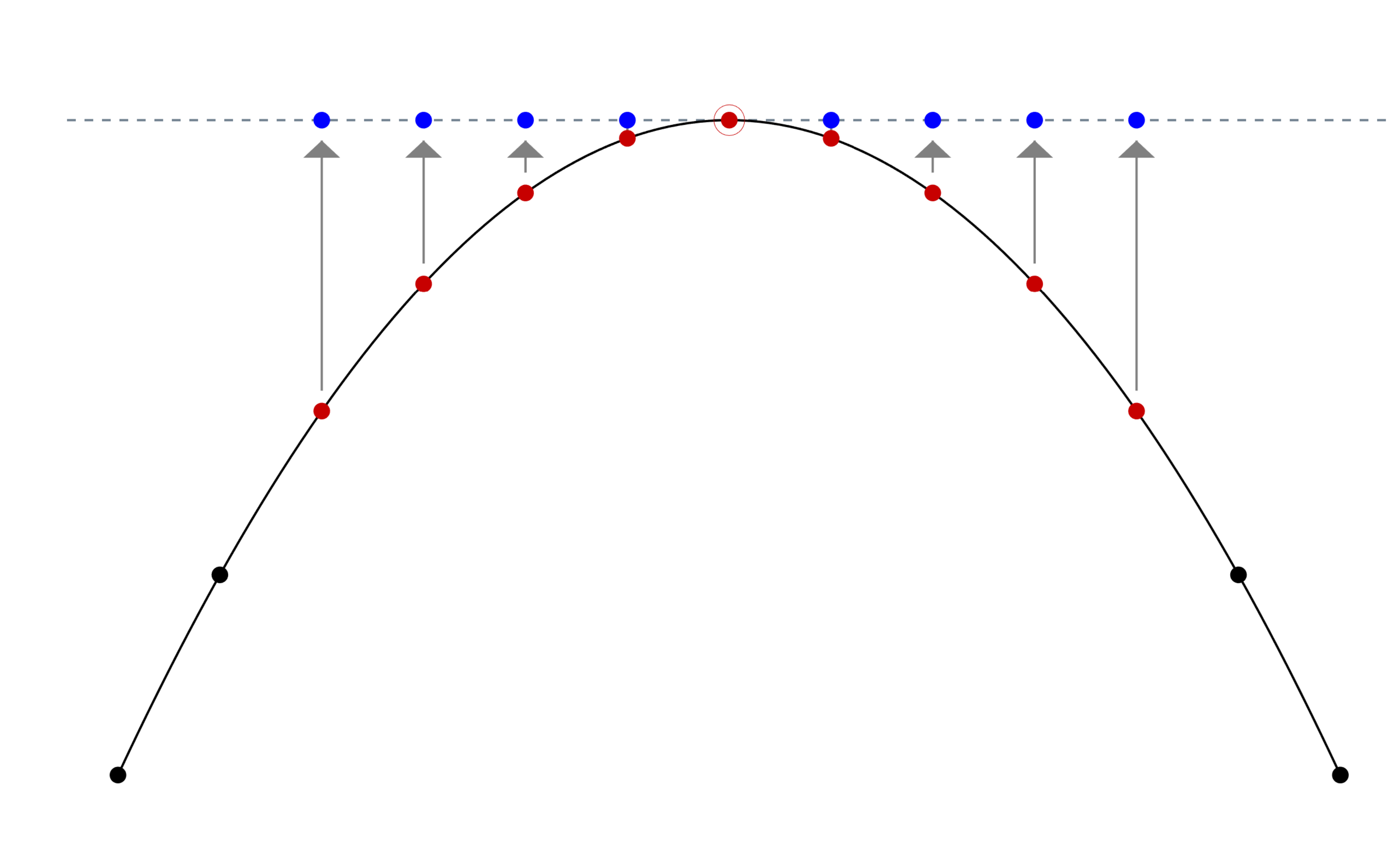}
  \caption{Projection to the tangent plane along the central normal. The central point is marked with an extra circle around it. Its neighbouring red points are projected to the tangent plane to the blue locations.}
  \label{Fig:TangentPlaneProjection}%
\end{figure}

For this, given a point $i$ on the surface, its neighbouring points $j\in S_i$ are projected to the tangent space $T_i$. The projection of the point $j\in S_i$ to $T_i$ will be referred to by $j_{T_i}$. This projection is done along the manifold normal $\vec{n}_i$. This process is illustrated for a $1$-dimensional manifold in $\mathbb{R}^2$ in Figure~\ref{Fig:TangentPlaneProjection}. The projection to the tangent space is not actually done at the numerical level. Only the distances between the central point $i$ and the projected locations $j_{T_i}$ are required. These are computed by rotating the distance vector $\delta\vec{x}_{ij} = \vec{x}_j - \vec{x}_i$. Using these distances, volumetric differential operators are computed in the tangent space. 

For a central point $i$, and projected locations $j_{T_i}$, usual $2$-dimensional gradients are computed along the $\vec{t}_{1,i}$ and $\vec{t}_{2,i}$ directions, given by
\begin{equation}
	\nabla_{T} \hat{u} \approx \widetilde{\nabla}_{T} \hat{u} = \left(\def\arraystretch{2.0} \begin{array}{c}
	\sum_{j\in S_i} c_{ij_T}^{t_1}\hat{u}_{j_T} \\													\sum_{j\in S_i}	c_{ij_T}^{t_2}\hat{u}_{j_T} \\			
		\end{array}\right)\,,
\end{equation}
where $\nabla_T$ denotes the gradient in the tangent plane, the $\tilde{\phantom{a}}$ overhead denotes the discrete derivative, and $\hat{u}_{j_{T}} = u_j$ denotes the function value at the projected location $j_{T_i}$. The stencil coefficients $c_{ij_T}$ are computed using a least squares approach while ensuring that monomials up to a certain order are differentiated exactly:
\begin{align}
	\sum_{j\in S_i}c_{ij_{T}}^{t_k}m_{j_{T}} &= \frac{\partial}{\partial t_k} m (\vec{x}_i) \qquad \forall m\in\mathcal{P}_{T}\,,\label{Eq:TP_Consistency}\\
	\text{min } J_i &= \sum_{j\in S_i} \left( \frac{c_{ij_{T}}^{t_k}}{W_{ij_{T}}} \right)^2\,, \label{Eq:TP_Min}
\end{align}
for $k=1,2$, where $\mathcal{P}_{T}$ is the set of monomials, usually up to order $2$, in $\vec{t}_{1,i}$ and $\vec{t}_{2,i}$ on the tangent plane. These are then rotated to obtain the numerical surface gradients. The surface derivatives are given by
\begin{equation}
	\label{Eq:SurfGrad_FullDefinition}
	\widetilde\nabla_{M,i} u = %
									\left(\def\arraystretch{2.0} \begin{array}{c}
									\sum_{j \in S_i}c_{ij}^{M,x} u_j \\												\sum_{j \in S_i}c_{ij}^{M,y} u_j \\		
									\sum_{j \in S_i}c_{ij}^{M,z} u_j \\		
										\end{array}\right)\,,
\end{equation}
where $c_{ij}^{M,x}$ are the stencil coefficients for the surface gradient in the $x$ direction, and similarly for the other directions. We have
\begin{equation}
	\left(\def\arraystretch{2.0} \begin{array}{c}
	c_{ij}^{M,x} \\														   			c_{ij}^{M,y} \\		
	c_{ij}^{M,z} \\		
		\end{array}\right) = %
	R^T 		\left(\def\arraystretch{2.0} \begin{array}{c}
	c_{ij}^{t_1} \\														   						c_{ij}^{t_2} \\		
	c_{ij}^{n} \\		
		\end{array}\right)\,,
\end{equation}
for $c_{ij}^n = 0$, and a rotation matrix $R$ given by
\begin{equation}
	R^T = \left(\begin{array}{ccc} 
	\vec{t}_1 & \vec{t}_2 & \vec{n} 
	\end{array}\right)\,,
\end{equation}
for column vectors $\vec{t}_1$, $\vec{t}_2$, and $\vec{n}$.
\subsubsection{Surface Divergence}~\\
We first rewrite the numerical surface gradients computed above as
\begin{equation}
	\label{Eq:SurfaceGradientConcise}
	\widetilde{\nabla}_{M,i}u = ( G_{1i} u, G_{2i} u, G_{3i} u)^T\,.
\end{equation}
Thus, $G_{1i}u = \sum_{j\in S_i}c_{ij}^{M,x}u_j$, and similarly for $G_{2i}$ and $G_{3i}$. Using this notation, the divergence of a vector valued function $\vec{v} = (v^1, v^2, v^3)$ defined on the manifold is given by
\begin{equation}
	\widetilde{\nabla}_{M}\cdot \vec{v} = \sum_{k=1,2,3} G_{ki} v^k\,.
\end{equation}
\subsubsection{Surface Laplacian}~\\
The numerical surface Laplacian  or Laplace Beltrami is given by
\begin{equation}
	\label{Eq:LapBeltNumericalDefinition}
	\widetilde{\Delta}_{M} u = \sum_{j\in S_i} c_{ij}^{\Delta_M}u_j\,.
\end{equation}
To determine the stencil coefficients $c_{ij}^{\Delta_M}$, we proceed in a manner similar to what was done above for the gradients. We first compute a $2$ dimensional volumetric Laplacian operator on the tangent plane:
\begin{align}
	\sum_{j\in S_i}c_{ij_{T}}^{\Delta_{T}}m_{j_{T}} &= \Delta_{T} m (\vec{x}_i) \qquad \forall m\in\mathcal{P}_{T}\,,\label{Eq:TPL_Consistency}\\
	\text{min } J_i &= \sum_{j\in S_i} \left( \frac{ c_{ij_{T}}^{\Delta_{T}} } {W_{ij_{T}} } \right)^2\,. \label{Eq:TPL_Min}
\end{align}
Since the Laplacian is rotationally invariant, the tangent plane Laplacian directly gives us the surface Laplacian $c_{ij}^{\Delta_M} = c_{ij}^{\Delta_T}\;$.

The verification of the concept and implementation of these differential operators for PDEs on stationary surfaces has been done in~\cite{Suchde2019}.

\subsection{Numerical Results}
All sparse linear systems arising from the discretization of the PDEs are solved using a BICGSTAB2 iterative solver \cite{BiCGSTAB}, without the use of any preconditioner. Monomials up to the order of $2$ are used in the computation of all the differential operators.

\subsubsection{Convection Diffusion Reaction Equation on an Expanding Sphere}~\\
We consider the case of an advection diffusion reaction equation on an evolving manifold
\begin{equation}
\label{Eq:AdvDiff}
	\frac{D\phi}{Dt} + \phi\nabla_M\cdot\vec{v} = \alpha\Delta_M\phi + f(\phi)\,,
\end{equation}
where $\vec{v}$ is the advection velocity of the surface, $\frac{D\phi}{Dt}$ is the advective surface material derivative, and $\alpha$ is the diffusion coefficient. Note that, in general, $\vec{v}$ can have a component both normal and tangential to the manifold. 

The time integration proceeds as follows. In each time step, first the Lagrangian movement step of Eq.\,\eqref{Eq:Move_o2} is performed. After that, integration of the remaining terms of Eq.\,\eqref{Eq:AdvDiff} is done with a first order implicit method
\begin{equation}
	\frac{\phi^{(n+1)}-\phi^{(n)}}{\Delta t}   + \phi^{(n+1)}\nabla\cdot\vec{v}^{\,(n+1)} = \alpha\Delta_M \phi^{(n+1)} + f(\phi^{(n+1)}) \,,
\end{equation}
where the bracketed superscripts indicate the time level. Discretization of the spatial operators is done as mentioned in Section~\ref{sec:DiffOp}. 

To check experimental orders of convergence of the proposed method, we consider the simple example of an expanding sphere. The initial unit sphere, discretized with irregularly spaced points, expands at a constant rate. The surface at time $t$ is given by
\begin{equation}
	x^2 + y^2 + z^2 = \left( r(t) \right)^2\,,
\end{equation}
with $r(t) = 1 + 0.5t$. The considered velocity field is $\vec{v}= 0.5\vec{n}$. The numerical solution is compared to a manufactured solution
\begin{equation}
	\phi_{\text{exact}}(\vec{x}, t) = \exp(-6t)xy\,.
\end{equation}
Since the exact solution and the location of the surface is known analytically, this serves as a good test case for the Lagrangian movement, as well as for the discretization of the spatial derivatives. The reaction term $f$ is determined analytically so that the manufactured solution satisfies Eq.\,\eqref{Eq:AdvDiff} with $\alpha=1$. This makes use of the following:
\begin{align}
	\nabla_M\cdot\vec{v} &= \frac{1}{r(t)}\,,\\
	\Delta_M\phi_{\text{exact}} &= -\frac{6\phi_{\text{exact} } }{ \left( r(t) \right) ^2}\,.
\end{align}
Relative errors in the numerical solution are measured at $t=1$ by
\begin{equation}
	\epsilon_2 = \left( \frac{ \sum_{i=1}^N \| \phi_i - \phi_{\text{exact}}(\vec{x}_i) \|^2}{ \sum_{i=1}^N \|\phi_{\text{exact}}(\vec{x}_i) \|^2} \right)^\frac12 \,.
\end{equation}
Note that the numerical errors found are not just a function of the numerical differential operators, but also of the Lagrangian movement and the addition/deletion of points to maintain point cloud regularity. 

The errors and convergence orders are tabulated in Table~\ref{Tab:ExpandingSphere_h} with a changing spatial discretization and a small time step of $\Delta t = 0.4h^2$. Convergence orders are measured with respect to the number of points at the initial time. A second order convergence is observed for the same. The numerical convergence with a changing time step is shown in Table~\ref{Tab:ExpandingSphere_dt} for $h=0.1$. In this case, first order convergence is observed.
\begin{table}[!htbp]
\caption{Errors at $t=1$ and convergence orders with $h$ for the expanding sphere test case. $h$ is the smoothing length, $N$ is the number of points in the entire domain at the initial time, $\epsilon_2$ is the relative error, and $r$ is the order of convergence of $\epsilon_{2}$.}
\centering
\label{Tab:ExpandingSphere_h}
	\begin{tabular}{|l|c|c|c|}
\hline
	$h$ & $N$ & $\epsilon_2$ & $r$ \\
\hline
&&& \\[\dimexpr-\normalbaselineskip+2pt]
	$0.4$  & $\phantom{11\,}480$  & $6.57\times 10^{-2}$ & $-$ \\
	$0.2$  & $\phantom{1}1\,806$ & $1.93\times 10^{-2}$ & $1.85$ \\
	$0.1$ & $\phantom{1}7\,446$ & $4.76\times 10^{-3}$ & $1.98$ \\
	$0.05$& $30\,054$ & $1.19\times 10^{-3}$ & $1.98$ \\
\hline	
	\end{tabular}
\end{table}
\begin{table}[!htbp]
\caption{Errors at $t=1$ and convergence orders with respect to the time step $\Delta t$ for the expanding sphere test case. $\epsilon_2$ is the relative error, and $r$ is the order of convergence of $\epsilon_{2}$.}
\centering
\label{Tab:ExpandingSphere_dt}
	\begin{tabular}{|l|c|c|}
\hline
	$\Delta t \times 10^{2}$ & $\epsilon_2$ & $r$ \\
\hline
&& \\[\dimexpr-\normalbaselineskip+2pt]
$1$   & $5.45\times 10^{-2}$ & $-$ \\
$1/2$  & $2.74\times 10^{-2}$ & $0.99$ \\
$1/4$  & $1.38\times 10^{-2}$ & $0.99$ \\
$1/8$  & $6.99\times 10^{-3}$ & $0.98$ \\
$1/16$ & $3.56\times 10^{-3}$ & $0.97$ \\
\hline	
	\end{tabular}
\end{table}

\subsubsection{Diffusion on a Deforming Half Pipe}
\begin{figure}
  \centering
  \includegraphics[width=0.49\textwidth]{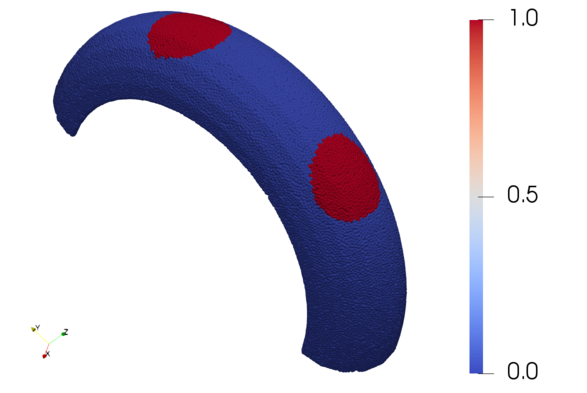}
  \includegraphics[width=0.49\textwidth]{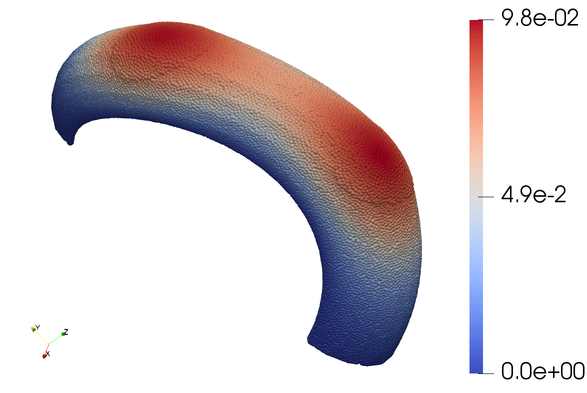}
  %
  \includegraphics[width=0.49\textwidth]{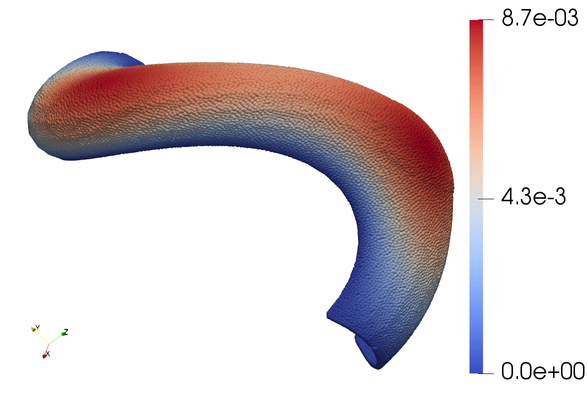}  
  \includegraphics[width=0.49\textwidth]{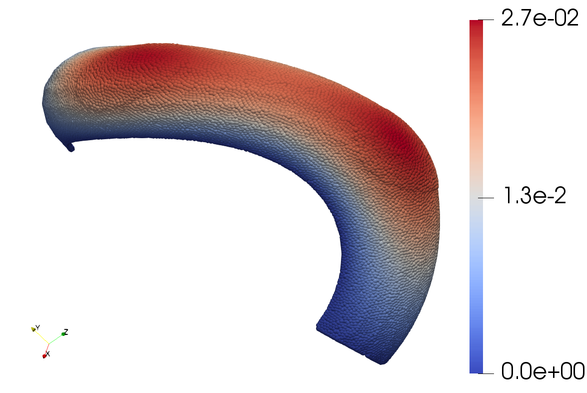}
  %
  \includegraphics[width=0.49\textwidth]{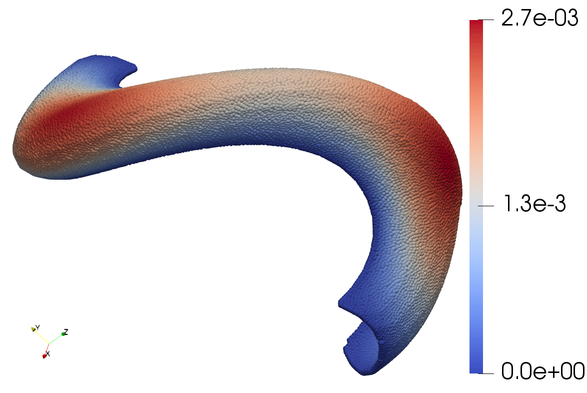}  
  \includegraphics[width=0.49\textwidth]{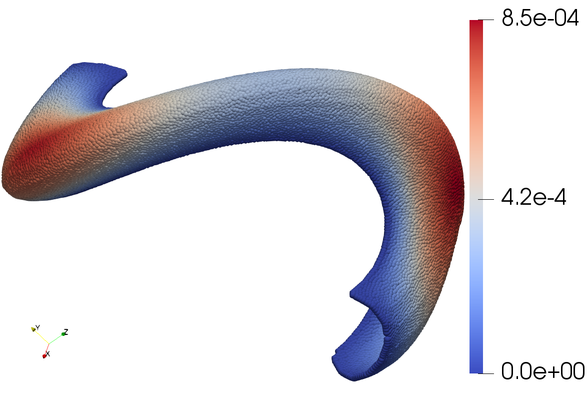}
  \caption{Diffusion on a deforming half pipe: Surface domain and solution shown at selected times. Clockwise from top left: At times $t = 0$~(top left), $t = 0.4$~(top right), $t = 0.8$~(middle right), $t = 1.2$~(middle left), $t = 1.6$~(bottom left) and $t = 2.0$~(bottom right). The colour represents the solution $\phi$.}
  \label{Fig:DeformingHalfPipe}%
\end{figure}

~\\
Next, we consider another convection diffusion reaction case according to Eq.\,\eqref{Eq:AdvDiff}. Unlike the previous case, here, the manifold has a boundary, the surface velocity $\vec{v}$ is dependent on the solution of the PDE, and also has a tangential component. Moreover, the deformation of the surface is more pronounced. Such cases of high surface deformation are important to consider here since these are the examples where the advantage of a moving meshfree method over a moving mesh method become evident. Under such high deformation, a moving mesh method would require a lot of remeshing, which would significantly slow down the simulation. The initial domain is taken to be one half of a pipe with radius $0.25$ bent in an arc of radius $0.75$, as shown in the first image in Figure~\ref{Fig:DeformingHalfPipe}. The velocity field is given by
\begin{equation}
	\label{Eq:HalfPipeVelocity}
	\vec{v} = 0.5\,\phi\, \vec{n} + \left( \begin{array}{c}
	-yz\\
	0\\
	xz\\
	\end{array} \right)\,,
\end{equation}
where  $\vec{n}$ is the unit normal vector, and $\phi$ is governed by Eq.\,\eqref{Eq:AdvDiff}. The first term in Eq.\,\eqref{Eq:HalfPipeVelocity} causes the surface to expand with the `temperature' $\phi$, while the second term causes an overall twisting motion with the two halves twisting in opposite directions. 
An initial condition with discontinuities is considered. At $t=0$, the domain has two `hot' regions with $\phi=1$, with the rest of the domain at $\phi=0$. We set $f= 0$ and $\alpha=0.2$ in Eq.\,\eqref{Eq:AdvDiff}. The domain is discretized with a smoothing length of $h=0.05$ that corresponds to an initial number of points of $N=11\,802$. Time integration is done as mentioned in the previous test case, with the velocity $\vec{v}$ taken as a function of $\phi$ at the previous time level. For a time step $\Delta t= 0.01$, Figure~\ref{Fig:DeformingHalfPipe} shows the evolution of the surface and the solution $\phi$. The initial half pipe twists in two different directions, with a slight expansion in the regions of high $\phi$, as $\phi$ diffuses. 


\subsubsection{Geometric Motion}
\begin{figure}
  \centering
  \includegraphics[width=0.45\textwidth]{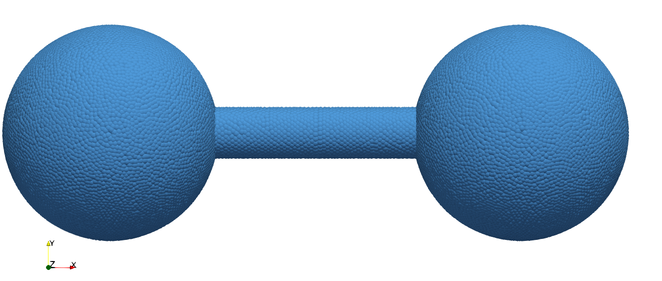}
  \includegraphics[width=0.45\textwidth]{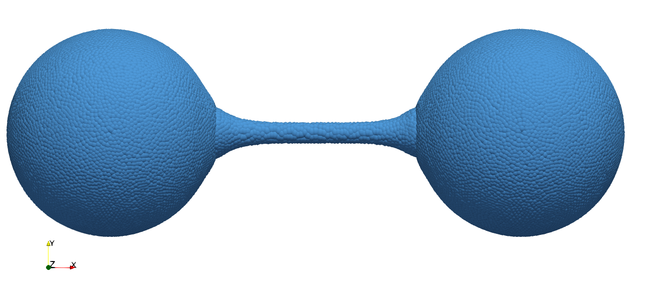}
  %
  \includegraphics[width=0.45\textwidth]{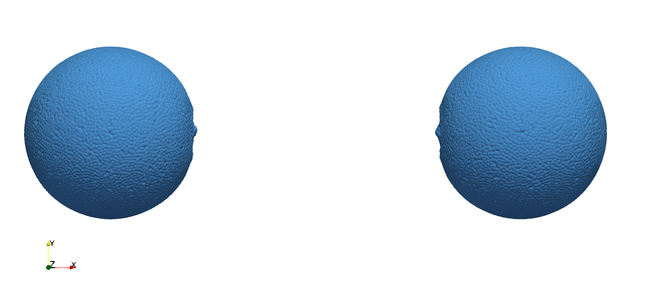}
  \includegraphics[width=0.45\textwidth]{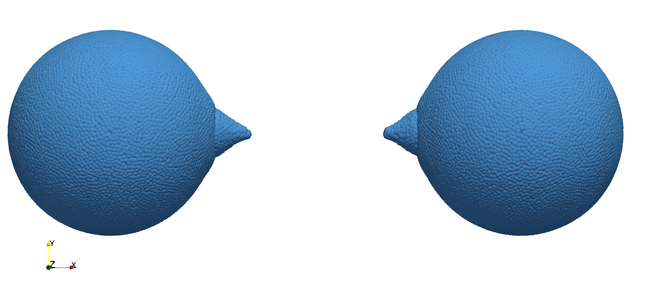}  
  \caption{Mean curvature flow for a dumbbell shape: Clockwise from top left: At times $t = 0$~(top left), $t = 0.05$~(top right), $t = 0.06$~(bottom right), and $t = 0.2$~(bottom left).}
  \label{Fig:MCF_Dumbbell}%
\end{figure}

\begin{figure}
  \centering
  \includegraphics[width=0.32\textwidth]{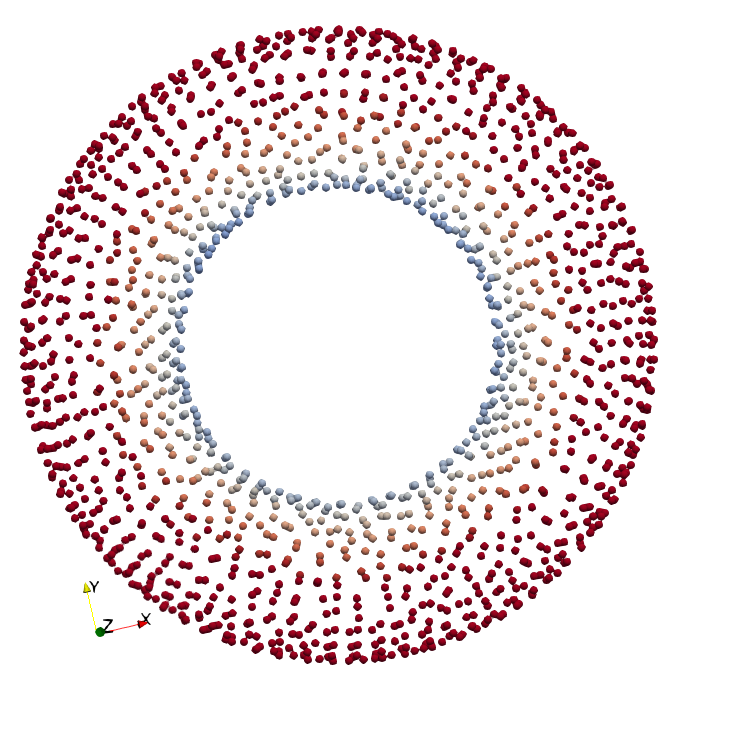}
  \includegraphics[width=0.32\textwidth]{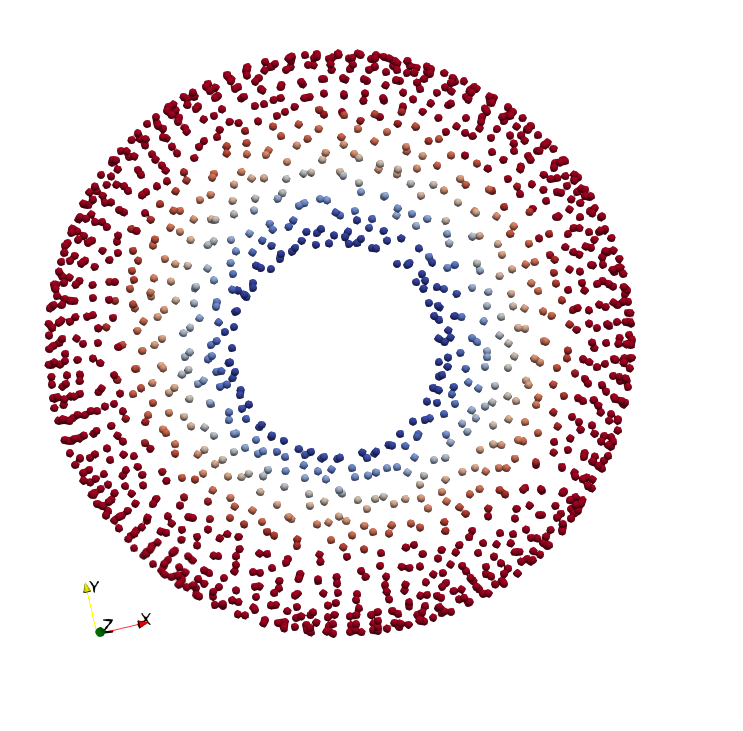}
  \includegraphics[width=0.32\textwidth]{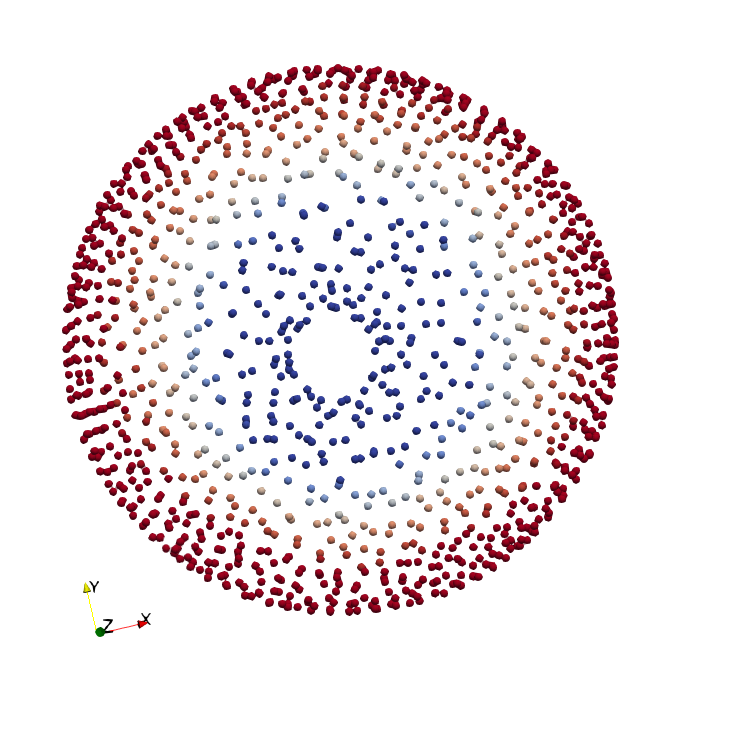}
  %
  \includegraphics[width=0.32\textwidth]{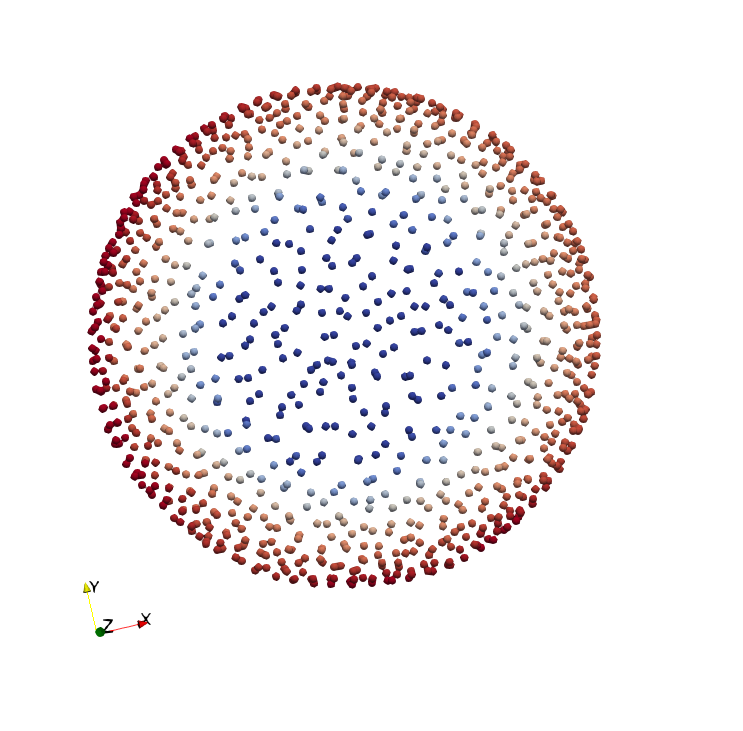}
  \includegraphics[width=0.32\textwidth]{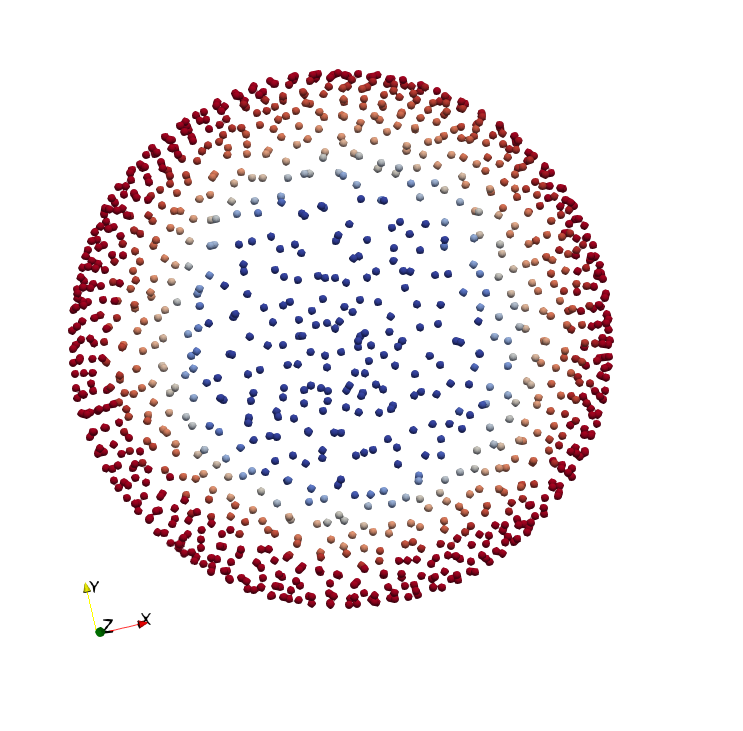}  
  \includegraphics[width=0.32\textwidth]{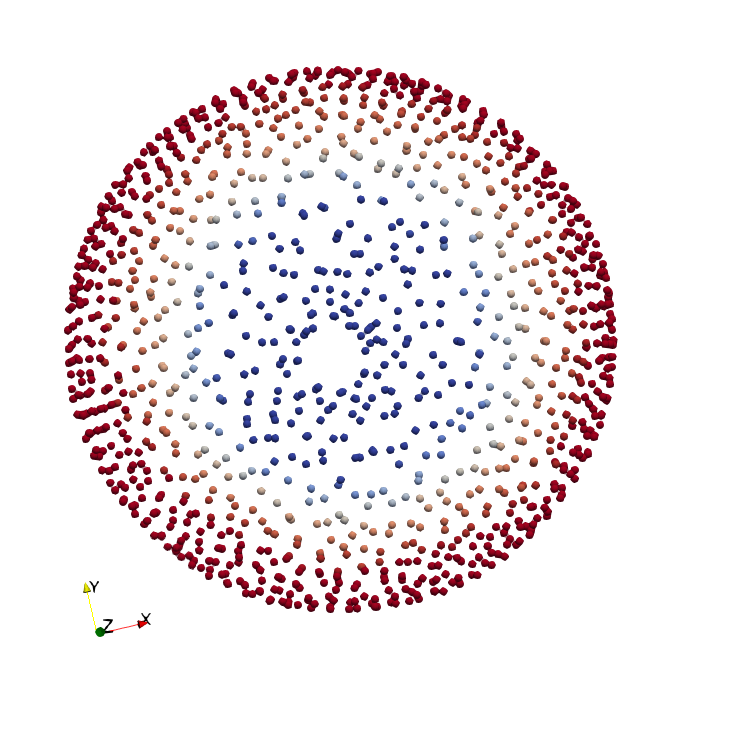}  
  \caption{Averaged mean curvature flow for a torus: Clockwise from top left: At times $t = 0$~(top left), $t = 4$~(top center), $t = 6.8$~(top right), $t = 7.05$~(bottom right), $t = 7.55$~(bottom center) and $t = 10$~(bottom left). The points are coloured according to the curvature $\kappa$ to make the results easier to visualize. The initial torus contracts and evolves to an ellipsoid. }
  \label{Fig:aMCF_Torus}%
\end{figure}

~\\
We consider motion of a surface dependent on its curvature. First, we consider the mean curvature flow~\cite{Colding2015}, in which the velocity of the surface is given by
\begin{equation}
	\label{Eq:MCF_vel}
	\vec{v} = \kappa \vec{n}\,,
\end{equation}
where $\vec{n}$ is the unit surface normal. The mean curvature $\kappa$ is given by
\begin{equation}
	\kappa = -\frac{1}{2}\nabla_M\cdot\vec{n} \,.
\end{equation} 
An initial dumbbell shape is used as the domain as shown in Figure~\ref{Fig:MCF_Dumbbell}. The dumbell is composed of two unit spheres with their centers $4$ units apart, connected by a cylindrical shape of radius $0.2$. There is a sharp change between the dumbbell ends and the handle. The simulation is done with $h=0.12$, which corresponds to an initial number of points of $N = 20\,015$, and a time step of $\Delta t = 0.005$. The results show the dumbbell shrinking, causing a neck pinch singularity, leading to the separation of the two ends, which evolve towards spheres as they contract. This matches the known behaviour of the dumbbell shape under mean curvature flow~\cite{Grayson1989, Mayer2000}.
As a second geometric motion test case, we consider the averaged mean curvature flow~\cite{Mayer2001}. The surface velocity is slightly modified from Eq.\,\eqref{Eq:MCF_vel} to
\begin{equation}
	\label{Eq:aMCF_vel}
	\vec{v} = (\kappa - \bar{\kappa}) \vec{n}\,,
\end{equation}
where $\bar{\kappa}$ denotes the global average of the curvature. The considered computational domain is the torus defined by
\begin{equation}
	\left( c - \sqrt{x^2 + y^2} \right)^2 + z^2 = a^2\,,
\end{equation}
for $c=3$ and $a = 1$. A coarsely discretized domain is used with $h=0.9$, which corresponds to $N=1\,490$ points at the initial state. The results for a time step of $\Delta t = 0.05$ are shown in Figure~\ref{Fig:aMCF_Torus}. The torus contracts towards the center, resulting in a topological change from genus $1$ to $0$, after which the domain slowly evolves towards a sphere. 

Both results show that the capability of the present method to capture topological change well. For both the dumbbell and the torus simulations, the contact handling algorithms introduced in Section~\ref{sec:Contact} detect the topological change, with the delete contact enforced, and the simulations do not need to be restarted after the singularities arise, which is often the case. 


An important point to note here is that to maintain stability, the values of $\kappa$ have to be smoothed before the velocity is computed. This requirement of smoothing has also been reported by other authors~\cite{Leung2011}. The smoothing is done with a wide Gaussian smoothing kernel on the support domain
\begin{equation}
	\tilde{\kappa}_i = \frac{\sum_{j\in S_i} \widetilde{W}_{ij} \kappa_j}{\sum_{j\in S_i} \widetilde{W}_{ij}} \,,
\end{equation}
with weights
\begin{equation}
	\widetilde{W}_{ij} = \exp\left(-\frac{\| \vec{x}_j - \vec{x}_i  \|^2}{h_i^2} \right)\,.
\end{equation}
%
%
%

\subsubsection{Wave Equation on an Evolving Surface}~\\
As the last experiment, we solve a wave equation on an evolving surface~\cite{Lubich2015}: 
\begin{equation}
	\label{Eq:WaveEqn}
	\frac{D}{Dt}\frac{D\phi}{Dt} + \frac{D\phi}{Dt}\nabla_M\cdot\vec{v} = c^2 \Delta_M \phi + f \,,
\end{equation}
where $\vec{v}$ is the surface velocity, and $f=f(\phi,\vec{x}, t)$ is a forcing term. Eq.\,\eqref{Eq:WaveEqn} is considered with initial conditions
\begin{align}
	\phi(\vec{x}, 0) &= g_1(\vec{x})\,,\\
	\frac{D\phi}{Dt}(\vec{x}, 0) &= g_2(\vec{x})\,. \label{Eq:WaveEqn_IC2}
\end{align}
Numerical integration in each time step begins with a movement step according to Eq.\,\eqref{Eq:Move_o2}, and is followed by
\begin{equation}
	\frac{ \phi^{(n+1)} - 2	\phi^{(n)} + 	\phi^{(n-1)} }{\left(\Delta t \right)^2} + \frac{	\phi^{(n+1)} - 	\phi^{(n)}}{\Delta t} \nabla_M \cdot \vec{v}^{\,(n+1)} = c^2 \Delta_M \phi^{(n+1)} + f^{(n)}\,.
\end{equation}
The numerical differential operators are computed as explained in Section~\ref{sec:DiffOp}, and $\phi^{(-1)}$ is computed based on Eq.\,\eqref{Eq:WaveEqn_IC2}. We start by validating the scheme with a manufactured solution
\begin{equation}
\label{Eq:Wave_MMS}
	\phi_{\text{exact}}(\vec{x}, t) =\sin(\omega t)xy\,.
\end{equation}
The considered surface at time $t$ is given by
\begin{equation}
	x^2 + y^2 + z^2 = \left( r(t) \right)^2\,,
\end{equation}
with $r(t) = 1 + 0.5t$. The considered velocity field is $\vec{v}= 0.5\vec{n}$. The forcing term and initial conditions are determined such that $\phi_{\text{exact}}$ in Eq.\,\eqref{Eq:Wave_MMS} satisfies Eq.\,\eqref{Eq:WaveEqn} with $c=1$. Similar to earlier, we note that the following holds:
\begin{align}
	\nabla_M\cdot\vec{v} &= \frac{1}{r(t)}\,,\\
	\Delta_M\phi_{\text{exact}} &= -\frac{6\phi_{\text{exact} } }{ \left( r(t) \right)^2}\,.
\end{align}
For a varying spatial discretization, convergence of the relative errors of the numerical solution as compared to the analytical solution are given in Table~\ref{Tab:WaveEqn}. The simulations use $\omega=\sqrt{6}$ in Eq.\,\eqref{Eq:Wave_MMS}, and a time step of $\Delta t = h / 20$. Convergence orders are measured with the number of points in the domain at initial time. The table shows that the numerical results match the manufactured solution. An experimental order of convergence of $2$ is observed which matches the expectation given the use of second order monomials in the computation of numerical derivatives.

\begin{table}[!htbp]
\caption{Errors at $t=1$ and convergence orders with $h$ for the wave equation test case with a manufactured solution. $h$ is the smoothing length, $N$ is the number of points in the entire domain at the initial time, $\epsilon_2$ is the relative error, and $r$ is the order of convergence of $\epsilon_{2}$.}
\centering
\label{Tab:WaveEqn}
	\begin{tabular}{|l|c|c|c|}
\hline
	$h$ & $N$ & $\epsilon_2$ & $r$ \\
\hline
&&& \\[\dimexpr-\normalbaselineskip+2pt]
	$0.4$  & $\phantom{11\,}480$  & $7.88\times 10^{-3}$ & $-$ \\
	$0.2$  & $\phantom{1}1\,806$ & $1.96\times 10^{-3}$ & $2.1$ \\
	$0.1$ & $\phantom{1}7\,446$ & $4.90\times 10^{-4}$ & $1.95$ \\
	$0.05$& $30\,054$ & $1.23\times 10^{-4}$ & $1.98$ \\
\hline	
	\end{tabular}
\end{table}

\begin{figure}
  \centering
  \includegraphics[width=0.32\textwidth]{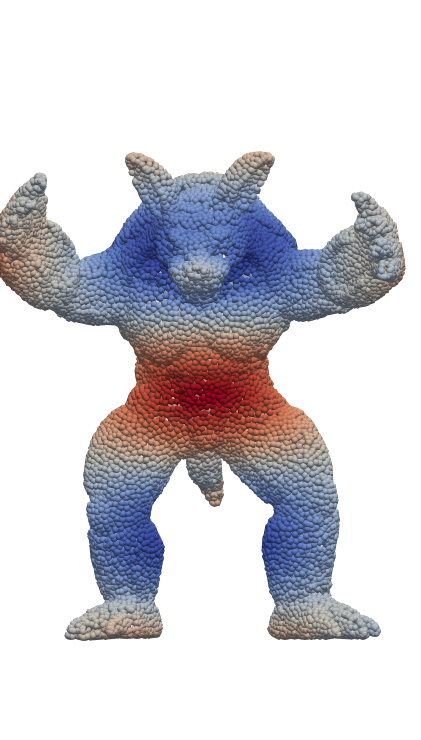}
  \includegraphics[width=0.32\textwidth]{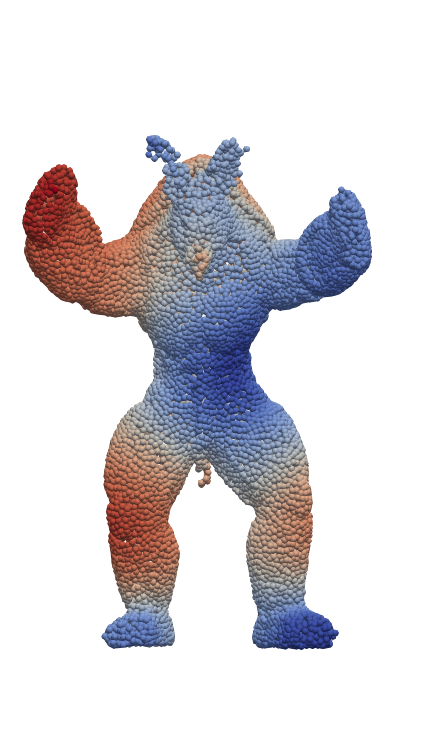}
  \includegraphics[width=0.32\textwidth]{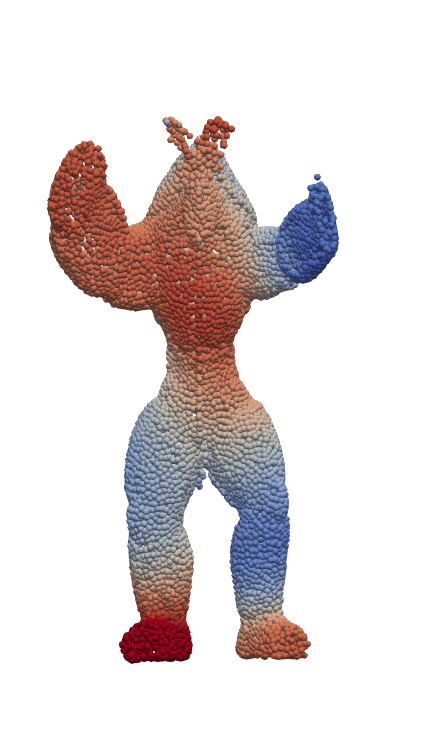}
  \includegraphics[width=0.32\textwidth]{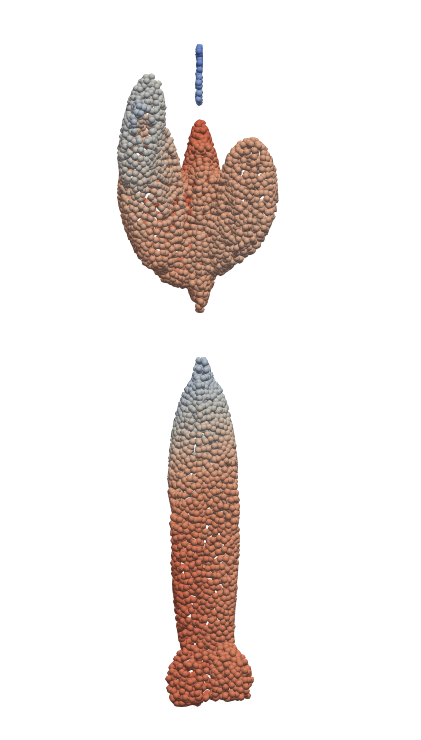}
  \includegraphics[width=0.32\textwidth]{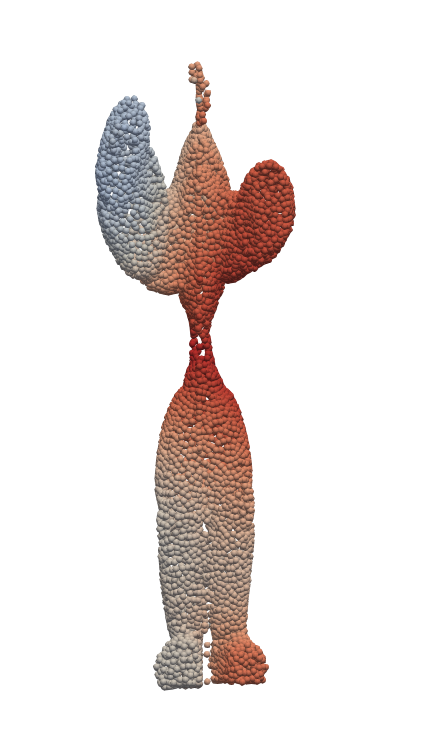}  
  \includegraphics[width=0.32\textwidth]{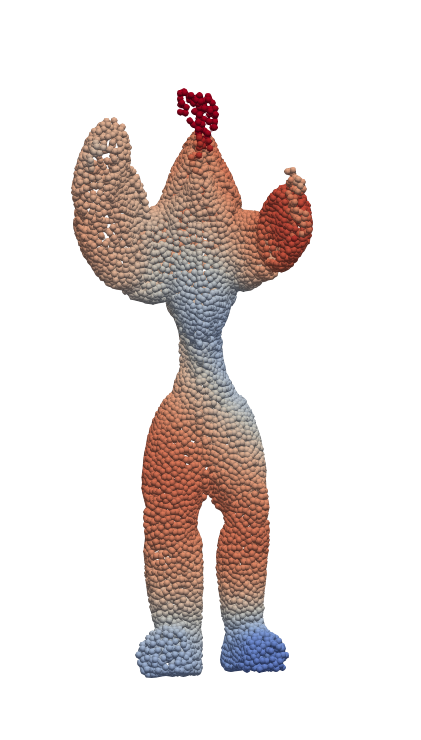}  
  \caption{Wave equation on a deforming armadillo: Clockwise from top left: At times $t = 0$~(top left), $t = 0.24$~(top center), $t = 0.48$~(top right), $t = 0.72$~(bottom right), $t = 0.96$~(bottom center) and $t = 1.2$~(bottom left). The points are coloured according to the solution $\phi$ of the wave equation. }
  \label{Fig:WaveArmadillo}%
\end{figure}

To further emphasize the capability of the present method to handle complex surfaces with large deformations, we consider the wave equation on a deforming Armadillo, with no forcing term. The considered initial conditions are
\begin{align}
	\phi^{(0)} &= -\cos(x)\sin(2y)\cos(z)\,,\\
	\phi^{(-1)} &=  -\cos(x-c\Delta t)\sin(2y-c\Delta t)\cos(z-c\Delta t)\,.
\end{align}
Furthermore, we set $f=0$, $c=7.0$, $\Delta t = 0.01$, and $h=0.2$ which corresponds to $N=13\,828$ points at the initial time. The velocity is given by
\begin{equation}
	\vec{v} = -0.3 \phi\vec{n} + \left( \begin{array}{c}
	-0.5x\\
	0.25y\\
	\\
	\end{array} \right) + \left( \begin{array}{c}
	-0.5 \text{sign}(x)\\
	0\\
	0\\
	\end{array} \right)\,.
\end{equation}
The simulation results are shown in Figure~\ref{Fig:WaveArmadillo}. The Armadillo stretches along its length as the wave travels in all directions. The collapsing of the Armadillo along its width triggers the contact detection algorithms, and results in a topological change in the surface. 

%
%
%
%
%
%

\section{Conclusion}
\label{sec:Conclusion}

We presented a novel fully Lagrangian framework to solve PDEs on moving manifolds. A surface is discretized with a point cloud, with no mesh or surrounding particles to discretize the volume around it. The method can handle arbitrary movements in the surface including those not prescribed a-priori and those that cause large deformations in the manifold. Distortions in the point cloud as a result of the evolution of the surface are fixed with purely local considerations. Similar to the volumetric case, the ease of fixing distortion is one of the biggest advantages of the moving surface point cloud over the moving surface mesh. We further presented a robust algorithm for contact handling for surface point clouds, that also enables simulations with topological changes in the surface.

To solve PDEs on evolving-in-time surfaces, this fully Lagrangian framework was used with a meshfree Generalized Finite Difference Method to approximate numerical derivatives. The Lagrangian framework was validated with several numerical experiments with different types of motion. And the applicability of the overall scheme was illustrated with applications to advection diffusion reaction equations on evolving surfaces, wave equations on evolving surfaces, and curvature-based geometric flow. This work forms a foundation for a purely meshfree Lagrangian method for PDEs on moving interfaces, with lots of potential applications ranging from shell mechanics to flow on curved surfaces.

A limitation of the present work is that differential operators are computed on a stationary point cloud at the beginning of each time step. As a result, achieving higher order accuracy for the overall scheme could be quite problematic. A possible solution to this could be to use more accurate operator splitting, as done, for example, in \cite{Leung2011}. Another improvement that needs to be explored is the correction of locations of points added to fill holes in the point cloud. In this paper, we considered the evolution of surfaces in both normal and tangential directions. However, evolution of particles within a (possibly stationary) surface was not considered. The extension to the present work to tackle the possibility of evolution within the manifold, for example for Lagrangian flow on a surface, will be considered in our future work. 


\section*{Acknowledgements}
The armadillo model is courtesy of the Stanford Computer Graphics Laboratory of Stanford University. 

%
%
%
%
%
%
%

\end{document}